\newif\ifieeetran
\IfFileExists{IEEEtran.cls}{%
  \ieeetrantrue
  \documentclass[10pt,draftclsnofoot,onecolumn]{IEEEtran}
}{%
  \ieeetranfalse
  \documentclass[10pt]{article}
  \usepackage[margin=25mm]{geometry}
}

\usepackage[T1]{fontenc}
\usepackage[utf8]{inputenc}
\usepackage{amsmath,amssymb,amsthm,mathtools,bm}
\usepackage{booktabs,tabularx,array,longtable}
\usepackage{graphicx}
\usepackage{float}
\usepackage{pdflscape}
\usepackage{enumitem}
\usepackage{setspace}
\usepackage{cite}
\usepackage{url}
\usepackage[hidelinks]{hyperref}
\usepackage[capitalise,noabbrev]{cleveref}
\ifieeetran\else
  \newenvironment{IEEEkeywords}{\par\noindent\textbf{Index Terms---}}{\par}
\fi

\graphicspath{{figures/}}
\allowdisplaybreaks[2]
\setlist{topsep=2pt,itemsep=0pt,parsep=0pt,partopsep=0pt}
\AtBeginDocument{%
  \setlength{\abovedisplayskip}{8pt plus 2pt minus 2pt}%
  \setlength{\belowdisplayskip}{8pt plus 2pt minus 2pt}%
  \setlength{\abovedisplayshortskip}{4pt plus 2pt minus 1pt}%
  \setlength{\belowdisplayshortskip}{6pt plus 2pt minus 1pt}%
  \setlength{\jot}{3pt}%
}
\theoremstyle{plain}

\newtheorem{definition}{Definition}
\newtheorem{remark}[definition]{Remark}
\newtheorem{implementationrule}[definition]{Implementation Rule}
\crefname{implementationrule}{implementation rule}{implementation rules}
\Crefname{implementationrule}{Implementation Rule}{Implementation Rules}

\usepackage{dsfont}
\usepackage{algorithm}
\usepackage{algorithmicx}
\usepackage{algpseudocode}
\makeatletter
\newenvironment{breakablealgorithm}
  {% \begin{breakablealgorithm}
   \begin{center}
     \refstepcounter{algorithm}% New algorithm
     \hrule height.8pt depth0pt \kern2pt% top rule
     \renewcommand{\caption}[2][\relax]{%
       {\raggedright\textbf{\ALG@name~\thealgorithm} ##2\par}%
       \ifx\relax##1\relax
         \addcontentsline{loa}{algorithm}{\protect\numberline{\thealgorithm}##2}%
       \else
         \addcontentsline{loa}{algorithm}{\protect\numberline{\thealgorithm}##1}%
       \fi
       \kern2pt\hrule\kern2pt
     }
  }{% \end{breakablealgorithm}
     \kern2pt\hrule\relax% bottom rule
   \end{center}
  }
\makeatother

\usepackage{cleveref}        % needed for \cref
\floatstyle{ruled}
\newcommand{\R}{\mathbb{R}}
\newcommand{\Nset}{\mathcal{N}}
\newcommand{\Eset}{\mathcal{E}}
\newcommand{\Iset}{\mathcal{I}}
\newcommand{\Bset}{\mathcal{B}}

\newcommand{\Gset}{\mathcal{G}}
\newcommand{\Tset}{\mathcal{T}}

\newcommand{\Mset}{\mathcal{M}}
\newcommand{\Nbr}[1]{\mathcal{N}_{#1}}
\newcommand{\Xset}{\mathcal{X}}
\newcommand{\DeltaID}{\Delta_{\mathrm{ID}}}
\newcommand{\DeltaDA}{\Delta_{\mathrm{DA}}}
\newcommand{\TUTC}{\mathbb{T}_{\mathrm{UTC}}}
\newcommand{\Tr}{^{\!\top}}
\newcommand{\defeq}{\vcentcolon=}
\newcommand{\pconv}{p^{\mathrm{conv}}}

\newcommand{\pconvDA}{p^{\mathrm{conv,DA}}}
\newcommand{\pcDA}{p^{\mathrm{c,DA}}}
\newcommand{\pdDA}{p^{\mathrm{d,DA}}}
\newcommand{\presDA}{p^{\mathrm{res,DA}}}
\newcommand{\pregDA}{p^{\mathrm{reg,DA}}}
\newcommand{\pgenDA}{p^{\mathrm{g,DA}}}
\newcommand{\ppumpDA}{p^{\mathrm{pump,DA}}}
\newcommand{\pplantDA}{p^{\mathrm{pl,DA}}}
\newcommand{\pflowDA}{p^{\mathrm{flow,DA}}}
\newcommand{\pconvSEC}{p^{\mathrm{conv,sec}}}
\newcommand{\pgenSEC}{p^{\mathrm{g,sec}}}
\newcommand{\ppumpSEC}{p^{\mathrm{pump,sec}}}
\newcommand{\pcstar}{p^{\mathrm{c},\star}}
\newcommand{\pdstar}{p^{\mathrm{d},\star}}
\newcommand{\pgenstar}{p^{\mathrm{g},\star}}
\newcommand{\ppumpstar}{p^{\mathrm{pump},\star}}
\newcommand{\pregstar}{p^{\mathrm{reg},\star}}
\newcommand{\Pmot}{P^{\mathrm{mot}}}
\newcommand{\Pav}{P^{\mathrm{av}}}
\newcommand{\Pres}{P^{\mathrm{res,max}}}

\newcommand{\etareg}{\eta^{\mathrm{reg}}}
\newcommand{\software}{\texttt{bahnstrom\_ems}}

\title{Distributed Intraday Energy Management for 16.7-Hz Railway Power Systems---Part II:\\
Software Realization}

\author{Navid Noroozi%
\thanks{N. Noroozi is with SIGNON Deutschland GmbH (DB InfraGO AG), Europaplatz 2, 10557 Berlin, Germany, \texttt{navid.n.noroozi@deutschebahn.com}}.}

\begin{document}
\maketitle

\begin{abstract}
\begin{singlespace}
This paper presents the executable Python implementation~\cite{noroozi2026bahnstromems} of the two-layer energy-management system developed mathematically in Part~I~\cite{norooziP1I}.
If needed, the software first converts hourly railway operation data into a 15-minute timetable and calculates the amount of gross motoring and available regenerative powers in a 15-minute resolution using the timetable handling layer.
Given the energy prices are available either cached or by sending request to ENTSO-E Transparency Platform through the price handling layer, the software constructs a deterministic hourly day-ahead plan with either in the so-called must-run converter commitment mode or in a sub-hourly feasibility certificate mode, where the latter demands more computational complexity.
Given the inputs from the timetable, day-ahead planning layers and initial measurements of the intraday-resolutions network's physical quantities, the intraday handling layer executes a risk-neutral intraday scenario model predictive controller by area-wise quadratic programming and consensus ADMM.
The retained S1 forecast layer combines a one-week seasonal-naive point forecast with causal residual-path resampling, while the current measured stage is identical in all scenarios.
We specify the day-ahead and intraday data contracts, absolute-UTC coupling, local-QP assembly, warm starts, residual tests, solver-status interpretation, physically
validated recovery paths, state and key-performance-index updates, and the railway-owned generation and pumping increment.
The Python package version, dependency lock, frozen-plan manifests, deterministic fixtures, run logs, and hash manifests make the software configuration executable and auditable.
The paper describes software behavior and numerical contracts; comparative performance and stress-test claims are reserved for Part~III~\cite{norooziP1III}.
\end{singlespace}
\end{abstract}

\begin{IEEEkeywords}
Railway power systems, energy management, reproducible software, Python programming, distributed convex optimization, model predictive control, point forecasting, unit commitment.
\end{IEEEkeywords}

\section{Introduction}
\label{sec:introduction}

The 16.7-Hz single-phase railway networks considered in Part~I~\cite{norooziP1I} contain geographically
distributed converter plants, railway-owned power generating and pumping units, wayside battery
energy storage systems (BESSs), renewable sources, and railway-power corridors.  Their energy
management is naturally divided between a slower day-ahead decision and a faster 15-minute
corrective dispatch.  Part~I formulates the corresponding risk-neutral scenario model
predictive control (MPC) problem, proves its convexity and sparse local structure, derives the
angle-copy consensus ADMM method, and establishes centralized equivalence under consensus
\cite{norooziP1I}.  Those mathematical statements deliberately abstract from such matters as
how an hourly public dataset is converted into quarter-hour inputs, how UTC indices are mapped,
which numerical result qualifies as convergence, or what the controller does after a local
solver fails.

The present Part~II closes that implementation gap for \software{} v0.9.1~\cite{noroozi2026bahnstromems}.  
Here show ``solver'' block, how the Part~I objects become a deterministic sequence of typed data, mixed-integer and quadratic programs, distributed coordination iterates, applied first-stage controls, physical-state updates, key performance indicators (KPIs), and persistent
run records.
The software consists of the following layers: (i) a causal and energy-conserving timetable synthesis
with a declared execution contract; (ii) a seasonal-naive, i.e. point, forecast and scenario interface;
(iii) live/deterministic-fixture economic inputs with identical optimization semantics;
(iv) a hourly day-ahead MILP and optional sub-hourly security block; 
(v) sparse area-QP assembly, consensus coordination, two levels of warm starting, and explicit failure
semantics;
(vi) a complete two-clock closed loop including first-stage application and KPI accumulation;
and (vii) a versioned package, plan-bank schema, result schema, and dependency lock.
The component, sequence, and class views are supplied as editable UML sources as well as rendered figures.

This paper is the middle part of a coordinated series.  Part~I is authoritative for the
physical model, objective, centralized scenario-MPC problem, decomposition, proofs, and
guarantee boundary~\cite{norooziP1I}.  The present Part~II is authoritative for the v0.9.1
software realization and execution semantics.
Part~III will illustrate the effectiveness of \software v0.9.1 through numerical results from two semi-synthetic benchmarks whose datasets given by~\cite{theiler2026}.

\section{Executable Scope, Clocks, and Software Contracts}
\label{sec:scope-contract}

We adopt the same notational and conceptual conventions as in Part~I~\cite{norooziP1I}.
For the sake of self-completeness, here we briefly recall the notation and the required concepts from~\cite{norooziP1I}. 

\subsection{Control areas and the two clocks}

Let the railway network be the connected graph
$\mathcal G=(\Nset,\Eset)$, where $i\in\Nset$ denotes an electrical control area and
$e\in\Eset$ a railway-power corridor.  The absolute time axis is
$\TUTC\subset\R_{\ge0}$ and all external timestamps are converted to timezone-aware UTC before
they enter a planning or control calculation.  The intraday sampling duration is
$\DeltaID=15\ \mathrm{min}=\tfrac14\ \mathrm h$,
$\DeltaDA=1\ \mathrm h$, $\nu\defeq\frac{\DeltaDA}{\DeltaID}=4$.
The closed-loop counter is $n\in\mathbb Z_{\ge0}$ and the $n$th control instant is the
absolute UTC time
$\tau_n: = \tau_0+n\DeltaID\in\TUTC$.
At $\tau_n$, the intraday horizon contains the relative prediction-stage indices $\Tset=\{0,\ldots,H-1\}$, $\tau_{n,t}=\tau_n+t\DeltaID$, $t\in\Tset$.

Thus $t$ is an integer stage index and never an absolute time or a sampling duration.
A stage is one 15-minute prediction interval together with the controls, states, and exogenous quantities assigned to that interval.  The day-ahead hour is the integer $h\in\mathcal H^{\mathrm{DA}}$ and is related to an intraday stage only through
\begin{equation}
h_n(t):=
\left\lfloor
\frac{\tau_{n,t}-\tau_0^{\mathrm{DA}}}{\DeltaDA}
\right\rfloor .
\label{eq:utc-map}
\end{equation}
Here $\tau_0^{\mathrm{DA}}$ is the absolute start of the active plan.  Neither $h$ nor $t$ is
used as a wall-clock timestamp.

The control instant $\tau_n$ is also the \emph{issue time}.  In plain language, it is the time
at which the controller takes its present information snapshot.  Mathematically, all current
measurements, forecast trajectories, asset-availability envelopes, and the selected active plan used
to form the problem at $\tau_n$ must be measurable with respect to the information set
$\mathcal F_{\tau_n}$.  The issue-time qualifier therefore identifies both the UTC origin and
the information boundary of one intraday calculation.

\subsection{Release and numerical scope}

The executable object described here is the Python package \software{} version 0.9.1~\cite{noroozi2026bahnstromems}.  It requires Python~3.11 or later.  NumPy and pandas hold numerical arrays
and UTC-indexed series; SciPy supplies sparse matrices and the HiGHS-backed \texttt{scipy.optimize.milp} interface; OSQP solves the convex intraday QPs and local ADMM subproblems~\cite{virtanen2020,huangfu2018,stellato2020}.
The package assembles the optimization matrices explicitly.
It does not require any external power-system simulation framework to construct the linearized railway-network model.

The release implements the following authority boundary:
\begin{itemize}[leftmargin=1.6em]
\item the day-ahead solver is centralized and replaceable, but its returned plan must satisfy
the declared \texttt{DayAheadPlan} contract;
\item the intraday controller is the risk-neutral finite-scenario controller of Part~I and is
coordinated through boundary-angle consensus;
\item S1 is the active transparent forecast construction described in this paper; 
\item railway-owned power plants have electrical generation and pumping envelopes; and
\item visualizations and exported files observe results after computation.
\end{itemize}

\subsection{Typed information flow}

The key software objects are summarized in \cref{tab:principal-contracts}.
Each object executes its associated data-class contract (container).

\begin{table}[htbp]
\centering
\caption{Key executable contracts.}
\label{tab:principal-contracts}
\small
\begin{tabularx}{\linewidth}{@{}p{3.4cm}X@{}}
\toprule
Object & Contract and role \\
\midrule
15-minute benchmark data
& One UTC-indexed base timetable containing $\Pmot_i$, $\Pav_i$, renewable availability,
prices, timetable provenance, and synthesis parameters.  The same object is supplied to
day-ahead security preparation and intraday execution. \\
\texttt{DayAheadPlan}
& Hourly UTC index, converter commitments and exchange schedule, BESS energy set-points,
peak targets, power plant schedules, and an optional quarter-hour security witness. \\
\texttt{ScenarioForecast}
& Forecast origin $\tau_n$, one \texttt{AreaExogenous} block per area, and a probability
vector $\pi$ with $\pi_m>0$ and $\sum_m\pi_m=1$. \\
\texttt{IntradayCoupling}
& Absolute-UTC mapping of an active plan into frozen commitment, optional gate-closed energy
equalities, BESS references, terminal floors, and peak targets for the current horizon. \\
\texttt{LocalQP}
& Sparse matrices $(P_i,q_i,A_i,l_i,u_i)$ and an explicit variable-index map for one area,
such that $\min\tfrac12x_i\Tr P_ix_i+q_i\Tr x_i$ subject to $l_i\le A_ix_i\le u_i$. \\
\texttt{AdmmResult}
& Area primal points, consensus and dual variables, penalty, residual history, convergence
flag, and classified local or outer failure information. \\
\texttt{ClosedLoopResult}
& Applied first stages, propagated BESS state, running peaks, realized KPIs, active-plan
history, scenario snapshots, solver diagnostics, and control-source labels. \\
\bottomrule
\end{tabularx}
\end{table}

\section{Input Preparation, Timetable Synthesis, S1 Forecasting, and Economics}
\label{sec:inputs}

This section explains how the software obtains the exogenous quantities whose availability is assumed by the optimization
model of Part~I.
It first gives the executable timetable-driven conversion from hourly targets to one common 15-minute resolution, then defines the
seasonal-naive point forecast and S1 scenario fan issued at each control instant, and finally
separates live price acquisition from deterministic economic fixtures.
The central objective is information consistency: day-ahead planning and intraday control must refer to the same identified base timetable and the same declared economic origin.

\subsection{One-time timetable-driven 15-minute synthesis}
\label{sec:synthesis-software}

For each day-ahead hour $h$, area $i$, and channel
$q\in\{\mathrm{mot},\mathrm{av}\}$, let $\overline P_i^q(h)$ be the hourly mean-power target.
The source dataset may additionally contain planned train counts by traffic category, train
mass or tonnage, and route context, but not the complete second-by-second positions of all
trains.  The \texttt{synth} package therefore combines representative train-leg calculations
with timetable templates.  Train resistance and electrical power follow the Davis-based
longitudinal model stated in Part~I~\cite{norooziP1I,davis1926,rochard2000}.

The synthesized raw clock-face shape is denoted
$\widetilde P_{i,h,k}^{q}\ge0$, where $k\in\{0,\ldots,\nu-1\}$ is the bin counter and $\nu$ captures the time-table resolution, e.g, $\nu = 4$ for quarter-hour.
In the complete-timetable case it is the binned integral of the positive motoring or braking channel.
In the public-data specialization it is the sum of planned category counts multiplied by representative acceleration, cruise, and braking
energies.
Arrival nodes place braking energy into the arrival bin, acceleration is placed in the following departure bin, cruise energy is distributed across the hour, and categories without declared clock-face nodes are uniform.
The normalized shape is
\begin{equation}
\widehat w_{i,h,k}^{q}=
\begin{cases}
\displaystyle
\widetilde P_{i,h,k}^{q}
\big/\sum_{j=0}^{\nu-1}\widetilde P_{i,h,j}^{q},
&\sum_j\widetilde P_{i,h,j}^{q}>0,\\[2mm]
1/\nu,&\text{otherwise}.
\end{cases}
\label{eq:normalized-template}
\end{equation}
The concentration parameter $c\in[0,1]$ forms the convex blend
$w_{i,h,k}^{q}(c)=\frac{1-c}{\nu}+c\widehat w_{i,h,k}^{q}$,
$P_{i,h,k}^{q}(c)=\nu\overline P_i^q(h)w_{i,h,k}^{q}(c)$.
This gives uniform quarter-hour values at $c=0$ and the full declared clock-face shape at $c=1$.
Nonnegativity and $\sum_k w_{i,h,k}^{q}=1$ imply the executable conservation assertion
\begin{equation}
P_{i,h,k}^{q}(c)\ge0,
\qquad
\frac1\nu\sum_{k=0}^{\nu-1}P_{i,h,k}^{q}(c)=\overline P_i^q(h),
\label{eq:executable-conservation}
\end{equation}
which is the property proved in Part~I.

\begin{algorithm}[H]
\caption{Timetable-driven construction of the common 15-minute resolution}
\label{alg:synthesis}
\begin{algorithmic}[1]
\Require hourly UTC index $\mathcal H$; targets
$\{\overline P_i^{q}(h)\}_{i\in\Nset,h\in\mathcal H,
q\in\{\mathrm{mot},\mathrm{av}\}}$; planned train counts and category map; representative
rolling-stock/route data; $\nu=4$; $c\in[0,1]$; synthesis revision $r_{\mathrm{syn}}$
\Ensure immutable 15-minute resolution
$\mathcal D_{15}=\{\Pmot_i(\tau),\Pav_i(\tau)\}_{i\in\Nset,\tau\in\mathcal T_{15}}$
and provenance $\mathcal P_{15}$
\State{$\Call{ValidateHourlyInputs}{\mathcal H,\overline P,c,\nu}\gets\mathtt{true}$}
\For{each declared traffic category $a$}
\State{$(E_a^{\mathrm{acc}},E_a^{\mathrm{cruise}},E_a^{\mathrm{brake}})
\gets\Call{TrainLegEnergy}{a,\text{rolling-stock data},\text{route data}}$}
\EndFor
\For{$i\in\Nset$, $h\in\mathcal H$, and
$q\in\{\mathrm{mot},\mathrm{av}\}$}
\State{$\widetilde P_{i,h,0:\nu-1}^{q}\gets
\Call{RawTemplate}{i,h,q,\text{counts},\text{category map},E^{\mathrm{acc}},
E^{\mathrm{cruise}},E^{\mathrm{brake}}}$}
\If{$\sum_{k=0}^{\nu-1}\widetilde P_{i,h,k}^{q}>0$}
\State{$\widehat w_{i,h,k}^{q}\gets
\widetilde P_{i,h,k}^{q}/\sum_{j=0}^{\nu-1}\widetilde P_{i,h,j}^{q}$,
$k=0,\ldots,\nu-1$}
\Else
\State{$\widehat w_{i,h,k}^{q}\gets1/\nu$, $k=0,\ldots,\nu-1$}
\EndIf
\For{$k=0,\ldots,\nu-1$}
\State{$w_{i,h,k}^{q}(c)\gets(1-c)/\nu+c\widehat w_{i,h,k}^{q}$}
\State{$P_{i,h,k}^{q}(c)\gets\nu\overline P_i^q(h)w_{i,h,k}^{q}(c)$}
\EndFor
\State{$\Call{Assert}{P_{i,h,k}^{q}(c)\ge0\ \forall k}$}
\State{$\Call{Assert}{\nu^{-1}\sum_{k=0}^{\nu-1}P_{i,h,k}^{q}(c)
=\overline P_i^q(h)}$}
\EndFor
\State{$\mathcal T_{15}\gets\Call{FlattenUTC}{\mathcal H,0:\nu-1}$}
\State{$\mathcal P_{15}\gets
\Call{Provenance}{\text{source IDs},\text{category parameters},c,\nu,r_{\mathrm{syn}}}$}
\Return{$(\mathcal D_{15},\mathcal P_{15})$}
\end{algorithmic}
\end{algorithm}

\begin{implementationrule}[Execute once, share twice]
\cref{alg:synthesis} is called during study preparation,
before either day-ahead planning or intraday dispatch.
The day-ahead input builder obtains its hourly means and its quarter-hour security profile from this same 15-minute object.
The
intraday controller later reads measurements and S1 history from the same object.
The
synthesis reissued only when the timetable, hourly targets, allocation, rolling-stock assumptions, or concentration value changes.
\label{rule:one-synthesis}
\end{implementationrule}

\subsection{Point forecasts and the S1 seasonal-naive predictor}
\label{sec:s1}

A point forecast is one predicted value of each uncertain quantity at each future stage, which is roughly speaking the center trajectory around which a finite scenario fan is formed.
Here we employ a specific, widely used point forecast approach called S1 seasonal forecaster~\cite{hyndman2021}.
Let $k_n$ be the integer location of $\tau_n$ in the stored 15-minute resolution and let $S=168\nu=672$ be the one-week lag in 15-minute samples.
The S1 point forecasts are seasonal persistence $
\widehat P_i^{\mathrm{mot}}(t\mid\tau_n) =P_i^{\mathrm{mot}}(k_n-S+t)$, $\widehat P_i^{\mathrm{av}}(t\mid\tau_n)=P_i^{\mathrm{av}}(k_n-S+t)$, $\widehat P_i^{\mathrm{res,max}}(t\mid\tau_n) =P_i^{\mathrm{res,max}}(k_n-S+t)$, for any stage $t\in\Tset$.
The exact residual construction and numerical guards below are project-specific engineering choices, not a claim of statistical calibration.

Import and export prices remain deterministic in S1.
All scenarios have equal probability
\begin{equation}
\pi_m=1/M, 
\label{eq:s1-probabilities}
\end{equation}
with $m\in\Mset=\{1,\ldots,M\}$.

%The provider interface is shared across benchmarks, while the residual transformation is a declared benchmark responsibility. 
% The SBB-3 15-minute realization uses one aggregate network residual path, realisation-normalized by the aggregate value, applies a common traction/regeneration multiplier clipped to $[0.7,1.3]$, and uses a separate milder photovoltaic perturbation.
%SBB-10 uses the area-wise construction \eqref{eq:sbb10-residual}--\eqref{eq:sbb10-res-scenario}.  Both retain the one-week centre, causal contiguous paths, equal probabilities, and the current-stage rule below.  This explicit distinction prevents a single shorthand formula from being falsely attributed to both executable backends.

After the provider returns, the closed-loop controller replaces stage zero in every scenario by the current measured values
\begin{align}
P_i^{\mathrm{mot}}(0,m)&\leftarrow P_{i,\mathrm{meas}}^{\mathrm{mot}}(\tau_n),\\
P_i^{\mathrm{av}}(0,m)&\leftarrow P_{i,\mathrm{meas}}^{\mathrm{av}}(\tau_n),\\
P_i^{\mathrm{res,max}}(0,m)&\leftarrow
P_{i,\mathrm{meas}}^{\mathrm{res,max}}(\tau_n),
\label{eq:stage-zero-overwrite}
\end{align}
for all $m \in \Mset$.
Current plant envelopes are likewise issued with a common stage zero.  Hence the action applied at $\tau_n$ is nonanticipative even though future stages branch.  Setting $\alpha_{\mathrm{res}}=0$ collapses the future traction and regenerative branches to the seasonal-naive centre, apart from the separately configured renewable perturbation if it is not also disabled.

%\begin{remark}[Scientific status of S1]
%S1 combines established seasonal persistence and contiguous residual-path resampling \cite{hyndman2021,kunsch1989,kaut2017}.  Its lag, normalization, common-factor choices, clipping bounds, and renewable perturbation are transparent benchmark definitions of this software.
%They are not estimated coverage probabilities and do not establish that the returned ensemble is probabilistically calibrated.  Part~III must evaluate any forecast-quality or decision-value claim against declared observations and comparators.
%\end{remark}

\subsection{Fixture and live economics}
\label{sec:economics}

Each converter plant exchanges energy with the public-grid bidding zone to which it is
connected.  From the zonal day-ahead price $\lambda_z(\vartheta)$, the software constructs
signed import and export prices
\begin{align}
c_i^{\mathrm{imp}}(\vartheta)
&=\lambda_{z(i)}(\vartheta)+a_i^{\mathrm{imp}},\\
c_i^{\mathrm{exp}}(\vartheta)
&=(1-\sigma_i^{\mathrm{exp}})\lambda_{z(i)}(\vartheta)-f_i^{\mathrm{exp}}.
\label{eq:price-transform}
\end{align}
The default engineering placeholders are an import adder of 12~EUR/MWh, an export haircut of
0.08, and an export fee of 2~EUR/MWh.  A pilot replaces them with the applicable contracts.
The export price is not clipped at zero: if it is negative, the objective term
$-c_i^{\mathrm{exp}}q_i^-$ correctly charges the controller for exporting.

The price handler has three ordered acquisition modes.  It first reads a local UTC cache.  If
coverage is missing and live acquisition is enabled, it requests ENTSO-E Transparency Platform
document A44 \cite{entsoe}, parses each period from its declared UTC start and resolution, and writes the
cache.  If a live token is intentionally unavailable and fixture operation is selected, it
constructs a deterministic seeded hourly series with declared monthly levels and a fixed
calendar shape.
The fixture is an offline test input, not observed hourly market evidence.
%SBB-10 applies a stricter publication boundary: its declared studies require complete cached Swiss prices and fail when that cache is absent or incomplete rather than silently substituting the fixture.

The full billing-period demand charge is assigned to the day-ahead peak variable.  The
intraday horizon uses only the configured pro-rata coefficient, because charging a short
receding horizon the full monthly peak tariff at every $\tau_n$ would repeatedly count one
billing obligation and would badly scale the local QPs.  Optional line-loss coefficients are
likewise part of the \texttt{EconomicConfig} provenance.  Changing live/fixture mode changes
the numerical coefficients and their evidential status, but not the optimization structure.

\section{Centralized Day-Ahead Planning Software}
\label{sec:day-ahead}

Here, we describe the replaceable slow layer that supplies the interface required by the
intraday controller.
It first states the deterministic hourly MILP implemented in v0.9.1,
then distinguishes the conservative must-run option from the sub-hourly-secured option, and
finally specifies the plan object and persistence checks.

\subsection{Hourly input and decision structure}

Let
\begin{equation}
\mathcal H^{\mathrm{DA}}=\{0,\ldots,H^{\mathrm{DA}}-1\}
\label{eq:da-hours}
\end{equation}
index the hourly plan beginning at $\tau_0^{\mathrm{DA}}$.  The deterministic input is one
hourly point trajectory for $\overline P_i^{\mathrm{mot}}(h)$,
$\overline P_i^{\mathrm{av}}(h)$, renewable availability, signed prices, and power plant
generation and pumping envelopes.  In the benchmark implementation these hourly railway
values are arithmetic means of the common 15-minute object produced by
\cref{alg:synthesis}; an external planner may instead receive an independently issued
day-ahead point forecast, provided its provenance is declared.

For a converter plant $i\in\Iset$, the principal hourly variables are commitment
$\delta_i^{\mathrm{DA}}(h)\in\{0,1\}$, start and stop indicators
$u_i^{\mathrm{DA}}(h),v_i^{\mathrm{DA}}(h)\in\{0,1\}$, net exchange
$\pconvDA_i(h)$, nonnegative import/export splits
$q_i^{+,\mathrm{DA}}(h),q_i^{-,\mathrm{DA}}(h)$, and peak target
$\bar p_i^{\mathrm{DA}}$.  They satisfy
\begin{align}
\pconvDA_i(h)
&=q_i^{+,\mathrm{DA}}(h)-q_i^{-,\mathrm{DA}}(h),
\label{eq:da-split}\\
\underline P_i^{\mathrm{conv}}\delta_i^{\mathrm{DA}}(h)
&\le \pconvDA_i(h)
\le \overline P_i^{\mathrm{conv}}\delta_i^{\mathrm{DA}}(h),
\label{eq:da-commitment-box}\\
\delta_i^{\mathrm{DA}}(h)-\delta_i^{\mathrm{DA}}(h-1)
&=u_i^{\mathrm{DA}}(h)-v_i^{\mathrm{DA}}(h),
\qquad
u_i^{\mathrm{DA}}(h)+v_i^{\mathrm{DA}}(h)\le1.
\label{eq:da-transition}
\end{align}
The first transition uses the declared pre-horizon state.  Sliding sums of $u_i^{\mathrm{DA}}$
and $v_i^{\mathrm{DA}}$ impose minimum on/off durations; a horizon sum of
$u_i^{\mathrm{DA}}$ limits starts; and adjacent values of $\pconvDA_i$ obey the
converter ramp limit.  These are the conventional unit-commitment constraints.  An import-only
converter plant has $\underline P_i^{\mathrm{conv}}=0$ and covers, at the dispatch level, an
aggregated rotary converter plant.  A reversible converter plant admits a negative lower
bound.

For each BESS $i\in\Bset$, the hourly charge, discharge, and end-of-hour energy variables
satisfy
$0\le \pcDA_i(h) \le \overline P_i^{\mathrm c}$, $0\le \pdDA_i(h) \le\overline P_i^{\mathrm d}$, 
$E_i^{\mathrm{DA}}(h)=E_i^{\mathrm{DA}}(h-1)
+\DeltaDA\!\left(\eta_i^{\mathrm c}\pcDA_i(h)
-\frac{\pdDA_i(h)}{\eta_i^{\mathrm d}}\right)$,
$\underline E_i\le E_i^{\mathrm{DA}}(h) \le\overline E_i$,
$E_i^{\mathrm{DA}}(H^{\mathrm{DA}}-1) \ge E_i^{\mathrm{end}}$.
An optional binary mode prevents simultaneous charge and discharge.  Renewable injection and
accepted regenerative power obey
$0\le\presDA_i(h)\le\overline P_i^{\mathrm{res}}(h)$, $0\le\pregDA_i(h)\le\overline P_i^{\mathrm{av}}(h)$.

For every railway-owned power plant $g\in\Gset_i$, the v0.9.1 MILP contains continuous
generation and pumping decisions
\begin{equation}
0\le\pgenDA_g(h)\le\overline P_g^{\mathrm g}(h),
\qquad
0\le\ppumpDA_g(h)\le\overline P_g^{\mathrm p}(h),
\qquad
\pplantDA_g(h)=\pgenDA_g(h)-\ppumpDA_g(h).
\label{eq:da-plant}
\end{equation}
The issued envelopes in \eqref{eq:da-plant} represent continuous electrical flexibility.
The plant quadratic regularization is retained in the MILP through a piecewise-linear epigraph.
For tangent points $\xi\in\Xi_g$ and $x_g(h)=\pplantDA_g(h)$, $z_g^{\mathrm{pl}}(h) \ge a_g\xi x_g(h)-\frac{a_g}{2}\xi^2$ for all $\xi\in\Xi_g$,
so minimizing $z_g^{\mathrm{pl}}$ gives the declared piecewise-linear approximation of
$a_gx_g^2/2$.  The default grid contains 17 tangents over the admissible net-power interval;
the exact quadratic is used in the intraday QP.

With the fixed network orientation and incidence coefficient $\sigma_{i,e}$, the centralized
hourly network constraints are
\begin{align}
\pflowDA_e(h) = & B_e\bigl(\theta_{f(e)}^{\mathrm{DA}}(h) -\theta_{s(e)}^{\mathrm{DA}}(h)\bigr),\nonumber\\
-\overline F_e \le & \pflowDA_e(h)\le\overline F_e,\nonumber\\
\pconvDA_i(h) +&\pdDA_i(h)-\pcDA_i(h) +\presDA_i(h)+\pregDA_i(h) +\sum_{g\in\Gset_i}\pplantDA_g(h) -\sum_{e\in\Eset}\sigma_{i,e}\pflowDA_e(h) = \overline P_i^{\mathrm{mot}}(h), \nonumber\\
\bar p_i^{\mathrm{DA}}&\ge q_i^{+,\mathrm{DA}}(h),\nonumber\\
\bar p_i^{\mathrm{DA}}&\ge\Pi_i^{\mathrm{prior}},
\label{eq:da-peak}
\end{align}
with absent-asset terms omitted and one reference angle fixed to zero.  Optional convex
piecewise-linear loss epigraphs may be added without changing the MILP class.

\subsection{Hourly economic objective}

The implemented deterministic objective can be written as
\begin{align}
J^{\mathrm{DA}}={}&
\sum_{h\in\mathcal H^{\mathrm{DA}}}\Bigg\{
\DeltaDA\sum_{i\in\Iset}
\left(c_i^{\mathrm{imp}}(h)q_i^{+,\mathrm{DA}}(h)
-c_i^{\mathrm{exp}}(h)q_i^{-,\mathrm{DA}}(h)\right)
\nonumber\\
&+\sum_{i\in\Iset}
\left(c_i^{\mathrm{su}}u_i^{\mathrm{DA}}(h)
+c_i^{\mathrm{nl}}\delta_i^{\mathrm{DA}}(h)\right)
+\DeltaDA\sum_{i\in\Bset}c_i^{\mathrm b}
\left(\pcDA_i(h)+\pdDA_i(h)\right)
\nonumber\\
&+\DeltaDA\sum_{g\in\Gset}
\left(c_g^{\mathrm g}\pgenDA_g(h)
+c_g^{\mathrm p}\ppumpDA_g(h)
+z_g^{\mathrm{pl}}(h)\right)
\nonumber\\
&+\DeltaDA\sum_i c_i^{\mathrm{cu}}
\left(\overline P_i^{\mathrm{res}}(h)-\presDA_i(h)\right)
+\DeltaDA\sum_i c_i^{\mathrm{rl}}
\left(\overline P_i^{\mathrm{av}}(h)-\pregDA_i(h)\right)
\Bigg\}
\nonumber\\
&+\sum_{i\in\Iset}c_i^{\mathrm{pk,period}}\bar p_i^{\mathrm{DA}}
+J_{\mathrm{loss}}^{\mathrm{DA}}.
\label{eq:da-objective}
\end{align}
The curtailment and regenerative-loss constants may be removed during matrix assembly and
added back when the reported objective is reconstructed.  This changes neither the optimizer
nor the physical dispatch.  The full billing-period peak coefficient appears once in
\eqref{eq:da-objective}, not once per intraday horizon.

\subsection{Must-run and sub-hourly-secured variants}
\label{sec:da-variants}

The software offers two distinct conservatism mechanisms.  In the hourly must-run variant,
selected converter plants satisfy
\begin{equation}
\delta_i^{\mathrm{DA}}(h)=1,
\qquad \forall h\in\mathcal H^{\mathrm{DA}}.
\label{eq:da-must-run}
\end{equation}
This fixes only converter-plant availability; it neither forces positive converter import nor
forces railway-owned power generation.  It is a simple commitment safeguard when no sub-hourly
certificate is requested.

The secured variant instead leaves converter commitment as a decision but appends a continuous
quarter-hour realization indexed by $(h,k)$,
$k\in\{0,\ldots,\nu-1\}$.  It repeats the exact linearized nodal balances, line-flow equations,
thermal limits, converter boxes, BESS dynamics and bounds, renewable and regenerative bounds,
and plant envelopes on the 15-minute grid supplied by \cref{alg:synthesis}.  The security
variables share the hourly integer commitments and peak target.  The principal linking
equalities are
\begin{align}
\frac1\nu\sum_{k=0}^{\nu-1}
\pconvSEC_i(h,k)
&=\pconvDA_i(h),
\label{eq:security-converter-link}\\
\frac1\nu\sum_{k=0}^{\nu-1}
\pgenSEC_g(h,k)
&=\pgenDA_g(h),
&
\frac1\nu\sum_{k=0}^{\nu-1}
\ppumpSEC_g(h,k)
&=\ppumpDA_g(h),
\label{eq:security-plant-link}\\
E_i^{\mathrm{sec}}(h,\nu-1)&=E_i^{\mathrm{DA}}(h),
&E_i^{\mathrm{sec}}(h,0)&\text{ continues from the preceding sub-step},
\label{eq:security-bess-link}\\
q_i^{+,\mathrm{sec}}(h,k)&\le\bar p_i^{\mathrm{DA}}.
\label{eq:security-peak-link}
\end{align}
The generation and pumping averages in \eqref{eq:security-plant-link} are constrained
separately; cancellation of simultaneous average generation and pumping cannot hide an
inconsistent hourly schedule.  No security binary is introduced.  The single set of hourly
converter binaries must admit both the hourly plan and the supplied quarter-hour realization.

The certificate proves feasibility for one declared base timetable and the initial states
used during planning.  It is not a robust certificate for every future forecast error, plant
outage, or BESS state reached after intraday recourse.  This limitation is why any security
sub-step proposed as an emergency first-stage action is revalidated against the current
measurement and state in \cref{sec:failure-semantics}.

\begin{remark}
The day-ahead planner is modular in a precise sense: an alternative optimization method may
replace \cref{alg:day-ahead} if it returns the same validated plan fields, units, UTC coverage,
and optional security witness.  The intraday controller neither imports the MILP variable map
nor assumes that HiGHS produced the plan.
\end{remark}

\begin{algorithm}[htbp]
\caption{Day-ahead plan construction}
\label{alg:day-ahead}
\begin{algorithmic}
\Require{network $\mathcal G$; planning anchor $\tau_0^{\mathrm{DA}}$; hourly point inputs
$\mathcal D^{\mathrm{DA}}$; initial converter/BESS states $x_0^{\mathrm{DA}}$; prior peak
$p_0^{\mathrm{pk}}$; economic configuration $\mathcal C^{\mathrm{econ}}$; planning mode
$\sigma_{\mathrm{DA}}\in\{\mathrm{must\mbox{-}run},\mathrm{secured}\}$; must-run map
$\delta^{\mathrm{MR}}$ or quarter-hour security profile $\mathcal D^{\mathrm{sec}}$}
\Ensure{validated UTC-indexed \texttt{DayAheadPlan} $\mathcal P^{\mathrm{DA}}$}
\State{$\Call{ValidateDAInputs}{\mathcal G,\tau_0^{\mathrm{DA}},\mathcal D^{\mathrm{DA}},
x_0^{\mathrm{DA}},p_0^{\mathrm{pk}},\mathcal C^{\mathrm{econ}}}\gets\mathtt{true}$}
\State{$\mathcal I^{\mathrm{DA}}\gets
\Call{AllocateHourlyVariables}{\mathcal G,\mathcal D^{\mathrm{DA}}}$}
\State{$J^{\mathrm{DA}}\gets\Call{AssembleObjective}{\eqref{eq:da-objective},
\mathcal I^{\mathrm{DA}},\mathcal C^{\mathrm{econ}}}$}
\State{$\mathcal C^{\mathrm{DA}}\gets
\Call{AssembleHourlyConstraints}{\eqref{eq:da-split}\text{--}\eqref{eq:da-peak},
\text{minimum-time/ramp/start limits},\mathcal I^{\mathrm{DA}}}$}
\State{$\mathcal I^{\mathrm{sec}}\gets\varnothing$}
\If{$\sigma_{\mathrm{DA}}=\mathrm{must\mbox{-}run}$}
\State{$\mathcal C^{\mathrm{DA}}\gets\mathcal C^{\mathrm{DA}}
\cup\Call{MustRunRows}{\delta^{\mathrm{MR}},\eqref{eq:da-must-run}}$}
\Else
\State{$\mathcal I^{\mathrm{sec}}\gets
\Call{AllocateSecurityVariables}{\mathcal G,\mathcal D^{\mathrm{sec}}}$}
\State{$\mathcal C^{\mathrm{sec}}\gets
\Call{AssembleSecurityRows}{\mathcal D^{\mathrm{sec}},
\eqref{eq:security-converter-link}\text{--}\eqref{eq:security-peak-link}}$}
\State{$\mathcal C^{\mathrm{DA}}\gets\mathcal C^{\mathrm{DA}}\cup\mathcal C^{\mathrm{sec}}$}
\EndIf
\State{$\mathcal M^{\mathrm{DA}}\gets
\Call{CanonicalMILP}{J^{\mathrm{DA}},\mathcal C^{\mathrm{DA}},
\mathcal I^{\mathrm{DA}},\mathcal I^{\mathrm{sec}}}$}
\State{$R^{\mathrm{DA}}\gets\Call{HiGHS}{\mathcal M^{\mathrm{DA}}}$}
\If{$R^{\mathrm{DA}}.\mathrm{status}\ne\mathrm{success}$}
\State{\textbf{raise} \texttt{DayAheadSolveError}$(R^{\mathrm{DA}}.\mathrm{status})$}
\EndIf
\State{$\mathcal P^{\mathrm{DA}}\gets
\Call{DecodePlan}{R^{\mathrm{DA}},\mathcal I^{\mathrm{DA}},\mathcal I^{\mathrm{sec}},
\tau_0^{\mathrm{DA}}}$}
\State{$\Call{ReconstructReportedCost}{\mathcal P^{\mathrm{DA}},
\text{constant curtailment/loss terms}}$}
\State{$\Call{ValidatePlan}{\mathcal P^{\mathrm{DA}},\mathcal G,
\mathcal D^{\mathrm{DA}},\mathcal D^{\mathrm{sec}}}\gets\mathtt{true}$}
\Return{$\mathcal P^{\mathrm{DA}}$}
\end{algorithmic}
\end{algorithm}

\subsection{Plan persistence and compatibility}

Plan preparation may be separated from closed-loop execution.  A frozen bank contains the
serialized plans and a compatibility manifest.  Version 0.9.1 compares the benchmark
identifier, stress case, concentration $c$, plant mode, planning hour, delivery duration,
look-ahead duration, time resolution, and secured/unsecured mode.
It also validates exact planning-anchor timestamps and horizon coverage.  

The day-ahead delivery duration is the nominal interval for which the plan schedules energy.
The look-ahead duration is an additional tail that gives an intraday horizon near the delivery
boundary valid plan data.
The closed-loop start time is an absolute control instant and need not equal the planning anchor.

\section{Intraday Scenario-MPC Software Realization}
\label{sec:intraday-software}

This section follows one issue-time problem from validated S1 trajectories to the collection
of sparse area QPs.  It first enforces the scenario and measurement contracts, then maps the
active day-ahead plan by absolute UTC, and finally explains the variable layout and physical
constraints assembled locally.  
We show that the executable matrices implement that formulation without an undeclared approximation.

\subsection{Issue-time exogenous-data validation}

For each area $i$, \texttt{AreaExogenous} contains arrays of shape $H\times M$ for
$P_i^{\mathrm{mot}}$, $P_i^{\mathrm{av}}$, optional $P_i^{\mathrm{res,max}}$, and each
plant-availability envelope.  Import/export prices have length $H$ and are deterministic under
S1.  Validation requires finite values, nonnegative physical availability, exact area and unit
identifiers, and the declared $H$ and $M$.  For every uncertain channel $v_i$,
\begin{equation}
v_i(0,1)=\cdots=v_i(0,M)
\label{eq:deterministic-first-stage}
\end{equation}
must hold after the measured overwrite \eqref{eq:stage-zero-overwrite}.  A shape or identifier
mismatch is an input error and is not delegated to OSQP.
The probability vector is validated separately:
\begin{equation}
\sum_{m\in\Mset}\pi_m=1, 
\label{eq:probability-contract}
\end{equation}
with $\pi_m>0$. The same $\pi_m$ multiplies every scenario-dependent local cost contribution.  Thus unequal
weights could be admitted by the optimizer without changing its structure, even though S1 in
this paper uses \eqref{eq:s1-probabilities}.

\subsection{Absolute-UTC day-ahead coupling}

The coupling builder evaluates \eqref{eq:utc-map} for every $t\in\Tset$.  It rejects a horizon
that precedes the plan, exceeds its UTC coverage, or maps to an absent area.  Four kinds of
information are then supplied to the local QP.

First, converter commitment is frozen:
$\underline P_i^{\mathrm{conv}}\delta_i^{\mathrm{DA}}(h_n(t))
\le\pconv_i(t,m)\le
\overline P_i^{\mathrm{conv}}\delta_i^{\mathrm{DA}}(h_n(t))$.
Second, if a day-ahead hour is explicitly gate-closed, let
$\mathcal T_{n,h}=\{t\in\Tset:h_n(t)=h\}$.  The covered-stage mean is constrained by
$\frac{1}{|\mathcal T_{n,h}|} \sum_{t\in\mathcal T_{n,h}}\pconv_i(t,m)=b_i^{\mathrm{DA}}(h)$ for all $m \in \Mset$.
The default has no firm hour, preserving intraday recourse.  Third, linear interpolation of
the hourly BESS set-points supplies a soft energy reference and the applicable terminal floor.
Fourth, $\bar p_i^{\mathrm{DA}}$ becomes the baseline above which the intraday peak slack is
charged.

The hourly plant schedules are mapped by the same UTC rule but are retained only as audit and
warm-start references.  They do not create an intraday equality, bound, or tracking penalty.
At $\tau_n$, current causal generation and pumping envelopes are authoritative.  This avoids
converting a day-ahead power plant schedule into a fictitious hard water or fuel
trajectory.

\subsection{Area variable ordering and sparse-QP assembly}

For each area $i$, \texttt{VariableIndex} assigns contiguous slices to the local decision
vector $x_i$.  The repeated stage--scenario blocks contain converter exchange and its
nonnegative split where present; BESS charge, discharge, and energy where present; renewable
injection; accepted regenerative power; peak slack; the area's voltage angle and each
neighbor-angle copy; and each local plant's generation and pumping.  Asset-free slices are not
allocated.  This design keeps the dimension proportional to the assets physically present in
the area. 
The assembler produces
\begin{equation}
\min_{x_i}\quad \frac12x_i\Tr P_i x_i+q_i\Tr x_i,
\qquad
l_i\le A_i x_i\le u_i.
\label{eq:local-qp}
\end{equation}
Its rows implement the Part~I local model:
\begin{itemize}[leftmargin=1.6em]
\item one nodal balance for every $(t,m)$, including motoring demand, accepted regenerative
power, converter exchange, BESS, renewable injection, power plant net injection, and all
incident local line flows;
\item BESS energy dynamics, physical bounds, issue-time initial energy, optional terminal
floor, and the soft day-ahead set-point term;
\item converter boxes, import/export split, ramp limits, frozen hourly commitment, optional
gate-closed mean equalities, and running-peak slack;
\item renewable and regenerative availability bounds;
\item generation and pumping bounds for each plant envelope, plus the exact convex quadratic
plant term;
\item local linearized line relations and thermal bounds using owned and copied boundary
angles; and
\item first-stage non-anticipativity equalities for every applied control variable.
\end{itemize}
The objective is the probability-weighted sum of market energy, BESS throughput, renewable
curtailment, regenerative spill, power plant cost, loss, BESS-reference, terminal, and
peak-above-target terms defined in Part~I.  The probability is applied during matrix assembly,
not afterward.  A small common Tikhonov curvature $\varepsilon_{\mathrm{reg}}=10^{-6}$ is added
to otherwise linear coordinates in both the distributed and centralized constructions; hence
it cannot create an artificial objective gap between them.

Variable ordering is time-major within each asset block, with scenario indices retained inside
each stage.  Constraint rows are likewise grouped by physical relation.  Because a BESS row
couples only adjacent energy stages, a nodal balance couples only variables at one $(t,m)$, and
a corridor relation couples only two angle coordinates, $A_i$ has horizon-local bandwidth and
sparsity.  The software constructs it by sparse triplets rather than by allocating a dense
matrix.

\begin{algorithm}[htbp]
\caption{Assembly of one issue-time scenario-MPC problem}
\label{alg:intraday-assembly}
\begin{algorithmic}
\Require{issue time $\tau_n$; active plan $\mathcal P_n^{\mathrm{DA}}$; persistent network
$\mathcal G$; measured BESS energy $E(\tau_n)$; running peak $p^{\mathrm{pk}}(\tau_n)$;
validated scenario forecast $\mathcal S_n$; issue-time plant envelopes
$\mathcal A_n^{\mathrm{pl}}$}
\Ensure{ordered local-QP map $\mathcal Q_n=\{Q_{i,n}\}_{i\in\Nset}$ and boundary-angle copy
topology $\mathcal T_n^{\theta}$}
\State{$\widehat{\mathcal S}_n\gets
\Call{MeasuredStageZero}{\mathcal S_n,\text{measurements at }\tau_n}$}
\State{$\Call{Assert}{\widehat{\mathcal S}_n\text{ satisfies }
\eqref{eq:deterministic-first-stage}\text{--}\eqref{eq:probability-contract}}$}
\State{$\mathcal G_n\gets\Call{CopyWithInitialBessState}{\mathcal G,E(\tau_n)}$}
\State{$\mathcal C_n^{\mathrm{ID}}\gets
\Call{MapPlanByUTC}{\mathcal P_n^{\mathrm{DA}},\tau_n,H,\eqref{eq:utc-map}}$}
\State{$\Call{ValidateCoverageAndCommitment}{\mathcal C_n^{\mathrm{ID}},
\mathcal P_n^{\mathrm{DA}}}\gets\mathtt{true}$}
\For{$i\in\Nset$}
\State{$\mathcal I_{i,n}\gets
\Call{AllocateVariableIndex}{i,H,M,\mathcal G_n}$}
\State{$P_{i,n},q_{i,n}\gets
\Call{AssembleLocalObjective}{i,\widehat{\mathcal S}_n,
\mathcal C_n^{\mathrm{ID}},\mathcal A_n^{\mathrm{pl}},p^{\mathrm{pk}}(\tau_n),
\mathcal I_{i,n}}$}
\State{$A_{i,n},l_{i,n},u_{i,n}\gets
\Call{AssembleLocalRows}{\text{physical},\text{non-anticipativity},
\text{day-ahead coupling},\text{loss cuts}}$}
\State{$Q_{i,n}\gets
\texttt{LocalQP}(P_{i,n},q_{i,n},A_{i,n},l_{i,n},u_{i,n},\mathcal I_{i,n})$}
\State{$\Call{ValidateLocalQP}{Q_{i,n}}\gets\mathtt{true}$}
\EndFor
\State{$\mathcal T_n^{\theta}\gets\Call{AngleCopyTopology}{\mathcal G_n,\mathcal Q_n}$}
\Return{$(\mathcal Q_n,\mathcal T_n^{\theta})$}
\end{algorithmic}
\end{algorithm}

\subsection{Centralized same-QP construction}

For validation, the software may stack the already assembled local QPs and add hard equality
rows between every boundary-angle copy and its global area angle.  It does not re-derive the
physical equations through a second model.  Consequently, the centralized reference and ADMM
solve the same regularized objective, local constraints, scenarios, probabilities, and
day-ahead coupling.  The only difference is whether angle agreement is imposed as one set of
hard rows or approached iteratively through consensus.  This construction is the executable
counterpart of the centralized-equivalence proof in Part~I.

\section{Consensus-ADMM Numerical Execution}
\label{sec:admm-software}

This section specifies the executable coordinator obtained from one ADMM iteration, a warm
start, convergence, and a local-solver failure.
It first presents the exact angle-consensus updates, then gives the outer stopping tests and adaptive-penalty rule, and finally separates
within-solve and cross-MPC warm starts from the failure and recovery semantics.
The purpose is to prevent a numerical status such as an iteration limit from being mistaken for a physical infeasibility certificate.

\subsection{Angle-copy topology and coordinator update}

For every area-angle owner $a\in\Nset$, the topology builder enumerates all local QPs that
contain a copy of $a$'s angle trajectory.  Let $S_i^{(a)}$ select that $H M$-element trajectory
from $x_i$ and let $z_a\in\R^{HM}$ be its global consensus value.  The implemented constraints
are $S_i^{(a)}x_i=z_a$, $i\in\{a\}\cup\Nbr{a}$.
Only boundary angles are exchanged.  Local converter, BESS, renewable, regenerative, plant,
and line-flow variables remain inside their owning area.

With unscaled multipliers $y_i^{(a)}$ and penalty $\rho^k>0$, area $i$ solves in parallel
\begin{align}
x_i^{k+1}=\arg\min_{x_i\in\Xset_i}\quad
&\frac12x_i\Tr P_i x_i+q_i\Tr x_i
\nonumber\\
&+\sum_{a\in\{i\}\cup\Nbr{i}}
\left[y_i^{(a),k\Tr}\left(S_i^{(a)}x_i-z_a^k\right)
+\frac{\rho^k}{2}\left\|S_i^{(a)}x_i-z_a^k\right\|_2^2\right].
\label{eq:admm-local-update}
\end{align}
The feasible set $\Xset_i=\{x_i:l_i\le A_ix_i\le u_i\}$ and the local physical objective do
not change across outer iterations.  For fixed $\rho^k$, only the augmented linear term changes.
The coordinator then computes
\begin{equation}
z_a^{k+1}=\frac{1}{d_a}
\sum_{i:\,S_i^{(a)}\ \mathrm{exists}}
\left(S_i^{(a)}x_i^{k+1}+\frac{y_i^{(a),k}}{\rho^k}\right),
\qquad
d_a=1+|\Nbr{a}|,
\label{eq:admm-z-update}
\end{equation}
fixes the reference-area trajectory to zero, and updates
\begin{equation}
y_i^{(a),k+1}=y_i^{(a),k}
+\rho^k\left(S_i^{(a)}x_i^{k+1}-z_a^{k+1}\right).
\label{eq:admm-dual-update}
\end{equation}
The reported objective is the sum of the original local objectives.  Augmented terms are not
included in the physical dispatch cost.

\subsection{Residuals, stopping tests, and penalty adaptation}

Stack all copy disagreements in
\begin{equation}
r^{k+1}=\operatorname{col}_{i,a}
\left(S_i^{(a)}x_i^{k+1}-z_a^{k+1}\right).
\label{eq:admm-primal-residual}
\end{equation}
The implementation computes the dual residual without constructing a dense consensus matrix:
\begin{equation}
\|s^{k+1}\|_2
=\rho^k\left(
\sum_{a\in\Nset}d_a\|z_a^{k+1}-z_a^k\|_2^2
\right)^{1/2}.
\label{eq:admm-dual-residual}
\end{equation}
Let $p=HM\sum_a d_a$ be the total number of angle-copy scalar equalities.  With
$\mathcal A_cx$ denoting the stacked local copies and $\mathcal B_cz$ their repeated global
values, the Boyd residual thresholds \cite{boyd2011} are
\begin{align}
\epsilon^{\mathrm{pri},k+1}
&=\sqrt p\,\epsilon^{\mathrm{abs}}
+\epsilon^{\mathrm{rel}}
\max\{\|\mathcal A_cx^{k+1}\|_2,\|\mathcal B_cz^{k+1}\|_2\},
\label{eq:eps-primal}\\
\epsilon^{\mathrm{dual},k+1}
&=\sqrt p\,\epsilon^{\mathrm{abs}}
+\epsilon^{\mathrm{rel}}\|\mathcal A_c\Tr y^{k+1}\|_2.
\label{eq:eps-dual}
\end{align}
In addition, the physical componentwise angle-copy discrepancy is
\begin{equation}
c_{\theta}^{k+1}=\|r^{k+1}\|_\infty.
\label{eq:max-consensus}
\end{equation}
An \texttt{AdmmResult} is marked converged only if
\begin{equation}
\|r^{k+1}\|_2\le\epsilon^{\mathrm{pri},k+1},
\qquad
\|s^{k+1}\|_2\le\epsilon^{\mathrm{dual},k+1},
\qquad
c_{\theta}^{k+1}<\epsilon_\infty^\theta,
\label{eq:admm-stop}
\end{equation}
when the physical gate $\epsilon_\infty^\theta$ is enabled.  Thus a small two-norm residual
cannot hide one unacceptable corridor-angle mismatch.

At configured intervals, residual balancing updates the bounded penalty:
\begin{equation}
\rho^{k+1}=\begin{cases}
\min\{\tau_\rho\rho^k,\rho_{\max}\},
&\|r^{k+1}\|_2>\mu_\rho\|s^{k+1}\|_2,\\
\max\{\rho^k/\tau_\rho,\rho_{\min}\},
&\|s^{k+1}\|_2>\mu_\rho\|r^{k+1}\|_2,\\
\rho^k,&\text{otherwise}.
\end{cases}
\label{eq:rho-update}
\end{equation}
The defaults are $\mu_\rho=10$ and $\tau_\rho=2$.  When $\rho$ changes, the local OSQP
objects are rebuilt because the augmented Hessians change; the current outer primal points are
then explicitly supplied as their initial points.  The unscaled multipliers in
\eqref{eq:admm-dual-update} are retained.

\begin{algorithm}[htbp]
\caption{Consensus ADMM with a classified local-QP contract}
\label{alg:admm}
\begin{algorithmic}
\Require{ordered local QPs $\mathcal Q=\{Q_i\}_{i\in\Nset}$; angle-copy topology
$\mathcal T^{\theta}$; options
$\mathcal O=(K_{\max},\epsilon^{\mathrm{abs}},\epsilon^{\mathrm{rel}},
\epsilon_\infty^\theta,\rho_{\min},\rho_{\max},\mu_\rho,\tau_\rho)$; optional warm state
$\mathcal W^0=(x^0,z^0,y^0,\rho^0)$}
\Ensure{\texttt{AdmmResult} $R^{\mathrm{ADMM}}$ with primal points, consensus/dual state,
residual history, and classified status}
\State{$(x^0,z^0,y^0,\rho^0)\gets\Call{InitializeADMM}{\mathcal Q,
\mathcal T^{\theta},\mathcal W^0,\mathcal O}$}
\For{$i\in\Nset$}
\State{$\mathcal O_i^{\mathrm{QP}}\gets
\Call{BuildOSQP}{Q_i,\rho^0,\text{unchanged local constraints}}$}
\State{$\Call{WarmStart}{\mathcal O_i^{\mathrm{QP}},x_i^0}$}
\EndFor
\For{$k=0,\ldots,K_{\max}-1$}
\For{$i\in\Nset$ in logically parallel order}
\State{$\Call{UpdateLinearTerm}{\mathcal O_i^{\mathrm{QP}},z^k,y_i^k,\rho^k}$}
\State{$(x_i^{k+1},\sigma_i,v_i^{\max})\gets
\Call{SolveOSQP}{\mathcal O_i^{\mathrm{QP}}}$}
\If{$\sigma_i=\mathrm{maximum\ iterations\ reached}$}
\State{$\mathcal O_i^{\mathrm{QP}}\gets\Call{ColdRebuild}{Q_i,z^k,y_i^k,\rho^k}$}
\State{$(x_i^{k+1},\sigma_i,v_i^{\max})\gets
\Call{SolveOSQP}{\mathcal O_i^{\mathrm{QP}}}$}
\EndIf
\If{$\sigma_i=\mathrm{primal\ infeasible}$}
\State{\textbf{raise} \texttt{AdmmLocalInfeasibility}$(i,k,\sigma_i)$}
\ElsIf{$\sigma_i\notin\{\mathrm{solved},\mathrm{solved\ inaccurate}\}$}
\State{\textbf{raise} \texttt{AdmmLocalSolveError}$(i,k,\sigma_i)$}
\ElsIf{$v_i^{\max}>2\times10^{-2}$}
\State{\textbf{raise} \texttt{AdmmLocalFeasibilityError}$(i,k,v_i^{\max})$}
\EndIf
\EndFor
\State{$z^{k+1}\gets\Call{ConsensusUpdate}{x^{k+1},y^k,\rho^k,
\mathcal T^{\theta},\eqref{eq:admm-z-update}}$}
\State{$z_{\mathrm{ref}}^{k+1}\gets0$}
\State{$y^{k+1}\gets
\Call{DualUpdate}{x^{k+1},z^{k+1},y^k,\rho^k,\eqref{eq:admm-dual-update}}$}
\State{$(r^{k+1},s^{k+1},c_{\theta}^{k+1})\gets
\Call{Residuals}{\eqref{eq:admm-primal-residual}\text{--}\eqref{eq:max-consensus}}$}
\State{$\mathcal H^{k+1}\gets
\Call{AppendHistory}{r^{k+1},s^{k+1},c_{\theta}^{k+1},J(x^{k+1}),\rho^k}$}
\If{the three tests in \eqref{eq:admm-stop} hold}
\Return{$\Call{AdmmResult}{x^{k+1},z^{k+1},y^{k+1},\rho^k,
\mathcal H^{k+1},\mathrm{converged}}$}
\EndIf
\If{$k+1$ is a penalty-update instant}
\State{$\rho^{k+1}\gets\Call{ResidualBalancing}{r^{k+1},s^{k+1},
\rho^k,\eqref{eq:rho-update}}$}
\If{$\rho^{k+1}\ne\rho^k$}
\State{$\{\mathcal O_i^{\mathrm{QP}}\}_{i\in\Nset}\gets
\Call{RebuildAndWarmStart}{\mathcal Q,\rho^{k+1},x^{k+1}}$}
\EndIf
\Else
\State{$\rho^{k+1}\gets\rho^k$}
\EndIf
\EndFor
\Return{$\Call{AdmmResult}{x^{K_{\max}},z^{K_{\max}},y^{K_{\max}},
\rho^{K_{\max}},\mathcal H^{K_{\max}},\mathrm{outer\_limit}}$}
\end{algorithmic}
\end{algorithm}

\subsection{Nested local-solver contract}
\label{sec:local-solver-contract}

OSQP solves each convex local QP \cite{stellato2020}.  A returned status of \emph{primal
infeasible} is kept distinct from \emph{maximum iterations reached}.  The former is a solver
certificate about the fixed local constraint set; the latter says only that the numerical
method did not meet its tolerances within its budget.  Since $\Xset_i$ is unchanged across
outer iterations, an iteration-limit status after earlier successful solves of the same local
constraints is not, by itself, evidence that the area's physical model became infeasible.

%The SBB-10 contract is shown in \cref{tab:sbb10-contract}.  Inner accuracy is matched to the outer physical purpose while an independent check prevents a nominal solver success from accepting an excessive unscaled row violation.

%\begin{table}[htbp]
%\centering
%\caption{Scale-appropriate SBB-10 local and outer solver contract.}
%\label{tab:sbb10-contract}
%\small
%\begin{tabular}{@{}lll@{}}
%\toprule
%Quantity & Value & Interpretation \\
%\midrule
%OSQP absolute tolerance & $10^{-3}$ & inner numerical stopping term \\
%OSQP relative tolerance & $10^{-4}$ & inner numerical stopping term \\
%OSQP maximum iterations & $100\,000$ & finite inner budget \\
%Scaled termination & disabled & stopping in unscaled row coordinates \\
%Maximum independent row violation & $2\times10^{-2}$ & reject if $l_i\le A_ix_i\le u_i$ fails beyond gate \\
%Cold retry & one & iteration-limit status only \\
%Outer $\epsilon^{\mathrm{abs}},\epsilon^{\mathrm{rel}}$ & $10^{-4},10^{-3}$ & Boyd tests \\
%Physical angle-copy gate & $10^{-2}$ rad & componentwise consensus \\
%Initial $\rho$; outer limit & $10$; $2000$ & SBB-10 runner configuration \\
%\bottomrule
%\end{tabular}
%\end{table}

The independent local violation is evaluated as $v_i^{\max}(x_i)=
\max_j\left\{[l_{i,j}-(A_ix_i)_j]_+, [(A_ix_i)_j-u_{i,j}]_+\right\}$.
A point with non-finite components or $v_i^{\max}>2\times10^{-2}$ is rejected even if OSQP reports a solved status.  If OSQP reaches its iteration limit, one optional cold retry constructs a fresh numerical workspace for the identical $(P_i,q_i,A_i,l_i,u_i)$ and the identical outer state.
It changes no constraint, scenario, tolerance, cost, or control authority.  Other failure statuses are not retried by this rule.

The exception classification is explicit.
The possible labels are
\path{primal_infeasible}, \path{dual_infeasible},
\path{numerical_iteration_limit}, \path{non_convex},
\path{numerical_accuracy_failure}, and \path{solver_failure}.  The record
contains the outer iteration, area identifier, original status, inner iteration count,
available primal/dual residuals, independent violation, and whether the cold retry occurred.

\subsection{Warm starts at two time scales}
\label{sec:warm-starts}

Two mechanisms are deliberately separated.  Within one call of \cref{alg:admm}, each area's
OSQP workspace reuses its previous local primal point while the augmented linear term changes.
This is an inner-solver acceleration and does not cross a control instant.

Across control instants, the complete final outer state $\mathcal W_n=(\{x_i\}_i,\{z_a\}_a,\{y_i^{(a)}\}_{i,a},\rho)$ is shifted by one stage.
For any time-major block
$v=(v(0),\ldots,v(H-1))$, define
\begin{equation}
\mathcal S_Hv=(v(1),\ldots,v(H-1),v(H-1)).
\label{eq:warm-shift}
\end{equation}
The shift is applied to every converter, BESS, renewable, regenerative, angle, line, loss, and
per-unit plant block, as well as every $z$ and $y$ trajectory; the final $\rho$ is retained.
The new problem's matrices and bounds remain authoritative, so the shifted point need not be
feasible.  Scenario labels are retained only as numerical positions: because S1 draws a new
fan at each $\tau_n$, scenario $m$ is not claimed to represent the same physical path across
adjacent horizons.

If dimensions, asset identifiers, or variable slices make the shift incompatible, the software
records the reason and starts the next ADMM solve cold.  A warm start is an acceleration aid and
cannot be allowed to turn an otherwise valid control problem into a control failure.

\subsection{Failure semantics and physical recovery}
\label{sec:failure-semantics}

An ADMM failure and a control failure are different events.  The first means that the
distributed path did not produce a point satisfying its numerical contract.  The second means
that no physically validated first-stage action is available.  The controller applies the
following ordered policy.

\begin{enumerate}[leftmargin=1.7em]
\item If ADMM satisfies \eqref{eq:admm-stop} and the configured outer physical limit, use its
first stage.  A routine centralized comparison, if scheduled, remains diagnostic and cannot
replace this successful distributed action.
\item After an ADMM failure, the controller may solve the centralized hard-consensus problem
formed from the identical local QPs.  A solution demonstrates feasibility of that numerical
model and separates distributed convergence from current-QP feasibility.  If the configured
recovery switch permits, its first stage is used as a rare centralized fallback.
\item If no centralized action is available and the active plan stores a security witness, the
controller extracts the quarter-hour sub-step for the current UTC instant.  It rechecks
converter commitment and exchange, BESS power and next energy, renewable/regenerative and
plant availability, network-wide and nodal balance, and every line limit.  A sub-step inside a
gate-closed hour is rejected when it cannot certify the remaining realized hourly-energy
equality after preceding intraday recourse.
\item If neither path supplies a valid action, strict mode terminates the study.  A non-strict
diagnostic mode may record non-finite control and KPI entries while holding dynamic state, but
it does not fabricate a dispatch or contaminate energy totals.
\end{enumerate}

The configuration may disable routine centralized comparisons, failure-triggered centralized
recovery, or both.  The persisted \path{control_source} therefore distinguishes the labels
\path{admm}, \path{centralized_fallback},
\path{secured_day_ahead_fallback}, and \path{no_feasible_fallback}.  Solver
diagnostics remain attached to the failed distributed attempt even when a later fallback
successfully supplies the applied action.

\section{Complete Two-Layer Closed-Loop Execution}
\label{sec:complete-loop}

This section places the preceding software components on one chronological path.  It first
shows the complete preparation-and-control algorithm, including the correct one-time call to
timetable synthesis, then explains the physical-state and KPI update after one applied stage,
and finally presents the S1 software execution workflow and call sequence.  The resulting
description is the executable counterpart of the receding-horizon policy in Part~I~\cite{norooziP1I}.

\subsection{Preparation and two-clock operation}

The planning clock fires at a configured daily UTC hour.  At the study start, the most recent
planning anchor not later than $\tau_0$ is selected.  A plan is then solved or loaded for every
planning anchor needed to cover the delivery interval and its MPC look-ahead tail.  The control
clock subsequently advances by $\tau_n=\tau_0+n\DeltaID\in\TUTC$.
Many control instants therefore share one active plan, and a plan refresh need not coincide with the first study instant.
Let $N_{\mathrm{ctl}}\in\mathbb Z_{\ge1}$ denote the declared number of 15-minute control instants in the study.

\begin{breakablealgorithm}
\caption{Complete v0.9.1 two-layer execution with S1, warm starts, and failure oracle}
\label{alg:complete-loop}
{\small
\begin{algorithmic}[1]
\Require{versioned study configuration $\mathcal C$; study interval
$[\tau_0,\tau_{N_{\mathrm{ctl}}})$; hourly railway data $\mathcal D_h$;
timetable/category data $\mathcal D^{\mathrm{tt}}$; economic source; network/assets;
$H$, $M$, $s_0$, $c$; planning, ADMM, centralized-oracle, fallback, and logging options}
\Ensure{frozen plan bank in preparation-only mode; otherwise validated
\texttt{ClosedLoopResult} $R^{\mathrm{CL}}$ and artifact set $\mathcal A^{\mathrm{out}}$}
\State{$\Call{ValidateConfiguration}{\mathcal C}\gets\mathtt{true}$}
\State{$\mathcal E\gets\Call{LoadEconomics}{\text{live/cache/fixture mode},\mathcal C}$}
\State{$(\mathcal D_{15},\mathcal P_{15})\gets
\Call{TimetableSynthesis}{\mathcal D_h,\mathcal D^{\mathrm{tt}},c}$
using \cref{alg:synthesis} exactly once}
\State{$\mathcal D_{15}\gets\Call{ApplyDeclaredStress}{\mathcal D_{15},\mathcal C}$}
\State{$(\mathcal G,\mathcal A^{\mathrm{pl}})\gets
\Call{BuildNetworkAndPlantAvailability}{\mathcal D_{15},\mathcal C}$}
\State{$\mathcal K^{\mathrm{DA}}\gets
\Call{RequiredPlanningAnchors}{\tau_0,N_{\mathrm{ctl}},H,\mathcal C}$}
\For{$\kappa\in\mathcal K^{\mathrm{DA}}$}
\If{day-ahead preparation is enabled}
\State{$(\mathcal D_\kappa^{\mathrm{DA}},\mathcal D_\kappa^{\mathrm{sec}})
\gets\Call{BuildDAInputsFromSameChronology}{\mathcal D_{15},\kappa,\mathcal C}$}
\State{$\mathcal P_\kappa^{\mathrm{DA}}\gets
\Call{DayAheadPlan}{\mathcal G,\mathcal D_\kappa^{\mathrm{DA}},
\mathcal D_\kappa^{\mathrm{sec}},\mathcal E,\mathcal C}$ using \cref{alg:day-ahead}}
\If{frozen-plan persistence is enabled}
\State{$\Call{PersistFrozenPlan}{\mathcal P_\kappa^{\mathrm{DA}},
\textproc{PlanManifest}(\mathcal C,\mathcal P_{15})}$}
\EndIf
\Else
\State{$(\mathcal P_\kappa^{\mathrm{DA}},\mathcal M_\kappa)
\gets\Call{LoadFrozenPlan}{\kappa,\mathcal C}$}
\State{$\Call{ValidateManifest}{\mathcal M_\kappa,
\text{benchmark},\text{stress},c,\text{plant mode},\text{clocks},
\text{horizons},\DeltaID,\text{security}}\gets\mathtt{true}$}
\EndIf
\EndFor
\If{prepare-day-ahead-plans-only is enabled}
\Return{$(\{\mathcal P_\kappa^{\mathrm{DA}}\}_{\kappa\in\mathcal K^{\mathrm{DA}}},
\text{plan-bank locations})$}
\EndIf
\State{$E(\tau_0),p^{\mathrm{pk}}(\tau_0),K_0,
\mathcal W\gets\Call{InitializeClosedLoop}{\mathcal G,\mathcal C}$,
with $\mathcal W=\varnothing$}
\For{$n=0,\ldots,N_{\mathrm{ctl}}-1$}
\State{$\tau_n\gets\tau_0+n\DeltaID$}
\State{$\mathcal P_n^{\mathrm{DA}}\gets
\Call{InstallCoveringPlan}{\tau_n,H,\{\mathcal P_\kappa^{\mathrm{DA}}\}}$}
\State{$y_n\gets\Call{Measure}{\Pmot,\Pav,\Pres,E,p^{\mathrm{pk}},\tau_n}$}
\State{$\mathcal A_n^{\mathrm{pl}}\gets
\Call{IssuePlantEnvelopes}{\mathcal A^{\mathrm{pl}},\mathcal F_{\tau_n}}$}
\State{$\mathcal S_n\gets\Call{S1}{\tau_n,H,M,s_0+n,\text{causal history}}$}
\State{$\mathcal S_n\gets
\Call{MeasuredStageZero}{\mathcal S_n,y_n,\mathcal A_n^{\mathrm{pl}}}$}
\State{$(\mathcal Q_n,\mathcal T_n^\theta)\gets
\Call{IntradayAssembly}{\tau_n,\mathcal P_n^{\mathrm{DA}},y_n,
\mathcal G,\mathcal S_n,\mathcal A_n^{\mathrm{pl}}}$ using
\cref{alg:intraday-assembly}}
\If{cross-MPC warm start is enabled and $\mathcal W$ is compatible}
\State{$\mathcal W_n^0\gets\mathcal S_H\mathcal W$ using \eqref{eq:warm-shift}}
\Else
\State{$\mathcal W_n^0\gets\varnothing$}
\EndIf
\State{$R_n^{\mathrm{cmp}},R_n^{\mathrm{cen}}\gets\varnothing$}
\State{$R_n^{\mathrm{ADMM}}\gets
\Call{ConsensusADMM}{\mathcal Q_n,\mathcal T_n^\theta,\mathcal W_n^0}$ using
\cref{alg:admm}}
\State{$\Call{PersistAttempt}{R_n^{\mathrm{ADMM}}}$; $\mathcal W\gets
\Call{RetainNumericallyUsefulState}{R_n^{\mathrm{ADMM}}}$}
\If{$R_n^{\mathrm{ADMM}}.\mathrm{converged}$ and its first stage is physically valid}
\State{$u_n\gets\Call{FirstStage}{R_n^{\mathrm{ADMM}}}$;
$\sigma_n^{u}\gets\mathtt{admm}$}
\If{$n$ is a scheduled centralized-comparison instant}
\State{$R_n^{\mathrm{cmp}}\gets\Call{CentralizedSameQP}{\mathcal Q_n}$
\textbf{ for validation only}}
\EndIf
\Else
\State{$R_n^{\mathrm{cen}}\gets
\Call{CentralizedRecoveryIfEnabled}{\mathcal Q_n,y_n}$}
\If{$R_n^{\mathrm{cen}}$ supplies a physically valid first stage}
\State{$u_n\gets\Call{FirstStage}{R_n^{\mathrm{cen}}}$;
$\sigma_n^{u}\gets\mathtt{centralized\_fallback}$}
\Else
\State{$u_n^{\mathrm{sec}}\gets
\Call{CurrentSecurityWitness}{\mathcal P_n^{\mathrm{DA}},\tau_n}$}
\If{$u_n^{\mathrm{sec}}$ is physically valid for $y_n$}
\State{$u_n\gets u_n^{\mathrm{sec}}$;
$\sigma_n^{u}\gets\mathtt{secured\_day\_ahead\_fallback}$}
\Else
\State{\textbf{raise} \texttt{NoPhysicallyValidatedFallback}$(n,\tau_n)$}
\EndIf
\EndIf
\EndIf
\State{$\Call{Apply}{u_n,\DeltaID}$}
\State{$E(\tau_{n+1})\gets\Call{BessUpdate}{E(\tau_n),u_n,\DeltaID}$}
\State{$p^{\mathrm{pk}}(\tau_{n+1})\gets
\Call{PeakUpdate}{p^{\mathrm{pk}}(\tau_n),u_n}$}
\State{$K_{n+1}\gets\Call{AccumulateRealizedKPI}{K_n,u_n,y_n,R_n^{\mathrm{ADMM}}}$}
\State{$\Call{RecordTick}{\tau_n,\mathcal P_n^{\mathrm{DA}},\mathcal S_n,u_n,
\sigma_n^u,R_n^{\mathrm{ADMM}},R_n^{\mathrm{cmp}},R_n^{\mathrm{cen}}}$}
\EndFor
\State{$R^{\mathrm{CL}}\gets\Call{BuildAndValidateClosedLoopResult}{
\{u_n\},\{E(\tau_n)\},\{p^{\mathrm{pk}}(\tau_n)\},K_{N_{\mathrm{ctl}}}}$}
\State{$\mathcal A^{\mathrm{out}}\gets
\Call{ExportReadOnlyArtifacts}{R^{\mathrm{CL}},\text{tables},\text{summary},
\text{schema},\text{hash manifest},\text{optional figures}}$}
\Return{$(R^{\mathrm{CL}},\mathcal A^{\mathrm{out}})$}
\end{algorithmic}
}
\end{breakablealgorithm}

The synthesis call near the beginning of \cref{alg:complete-loop} is deliberately outside
the control loop.  It prepares one identified base trajectory before day-ahead planning.
S1 is called inside the control loop because its causal path pool and origin change with
$\tau_n$.  Confusing these two operations would
either regenerate a benchmark during control or freeze a forecast that should be reissued.

\subsection{Applied-state propagation and KPI accumulation}

Let the selected first-stage values at $\tau_n$ carry a superscript $\star$.  Only BESS energy
is a propagated plant state in v0.9.1:
\begin{equation}
E_i(\tau_{n+1})=E_i(\tau_n)
+\DeltaID\left(
\eta_i^{\mathrm c}\pcstar_i(0)
-\frac{\pdstar_i(0)}{\eta_i^{\mathrm d}}
\right).
\label{eq:realized-bess-update}
\end{equation}
The update is checked against $[\underline E_i,\overline E_i]$ before the next problem is
formed.  The metered import peak evolves as
\begin{equation}
\Pi_i(\tau_{n+1})=
\max\{\Pi_i(\tau_n),q_i^{+,\star}(0)\}.
\label{eq:running-peak}
\end{equation}
Hence a realized peak is not forgotten when the horizon recedes.

For example, realized market and power plant costs accumulate as
\begin{align}
J_{\mathrm{market}}(N)
&=\sum_{n=0}^{N-1}\DeltaID\sum_{i\in\Iset}
\left(c_i^{\mathrm{imp}}(\tau_n)q_i^{+,\star}(0)
-c_i^{\mathrm{exp}}(\tau_n)q_i^{-,\star}(0)\right),
\label{eq:realized-market}\\
J_{\mathrm{plant}}(N)
&=\sum_{n=0}^{N-1}\DeltaID\sum_{g\in\Gset}
\left[c_g^{\mathrm g}\pgenstar_g(0)
+c_g^{\mathrm p}\ppumpstar_g(0)
+\frac{a_g}{2}\left(\pgenstar_g(0)-\ppumpstar_g(0)\right)^2\right].
\label{eq:realized-plant}
\end{align}
Regenerative spill and the recovery ratio are
$E^{\mathrm{spill}}(N)
=\sum_{n=0}^{N-1}\DeltaID\sum_i
\left(\Pav_i(\tau_n)-\pregstar_i(0)\right)$,
$\etareg(N)
=\frac{\sum_n\DeltaID\sum_i\pregstar_i(0)}
{\sum_n\DeltaID\sum_i\Pav_i(\tau_n)}$, 
when the denominator is positive.
Plant telemetry separately reports generation energy, pumping energy, net electrical energy, opportunity cost, and envelope margins.
It reports no reservoir-energy KPI because no such state exists.

Failure-aware numerical KPIs retain unsuccessful attempts.
If $N_{\mathrm{att}}$ ADMM solves were attempted and $N_{\mathrm{succ}}$ converged,
$r_{\mathrm{ADMM}}=\frac{N_{\mathrm{succ}}}{N_{\mathrm{att}}}$,
$\bar k_{\mathrm{all}}=\frac{1}{N_{\mathrm{att}}}\sum_{n\in\mathrm{att}}k_n$,
$\bar k_{\mathrm{succ}}=\frac{1}{N_{\mathrm{succ}}}\sum_{n\in\mathrm{succ}}k_n$.
An outer-limit failure remains in $r_{\mathrm{ADMM}}$ and $\bar k_{\mathrm{all}}$; it is
excluded only from the statistic explicitly named successful.
The separate \texttt{solve\_failed} and \texttt{control\_failed} fields implement the distinction stated
in \cref{sec:failure-semantics}.

\subsection{Software execution workflow}

\begin{figure}[htbp]
\centering
\includegraphics[height=0.7\textheight,keepaspectratio]
{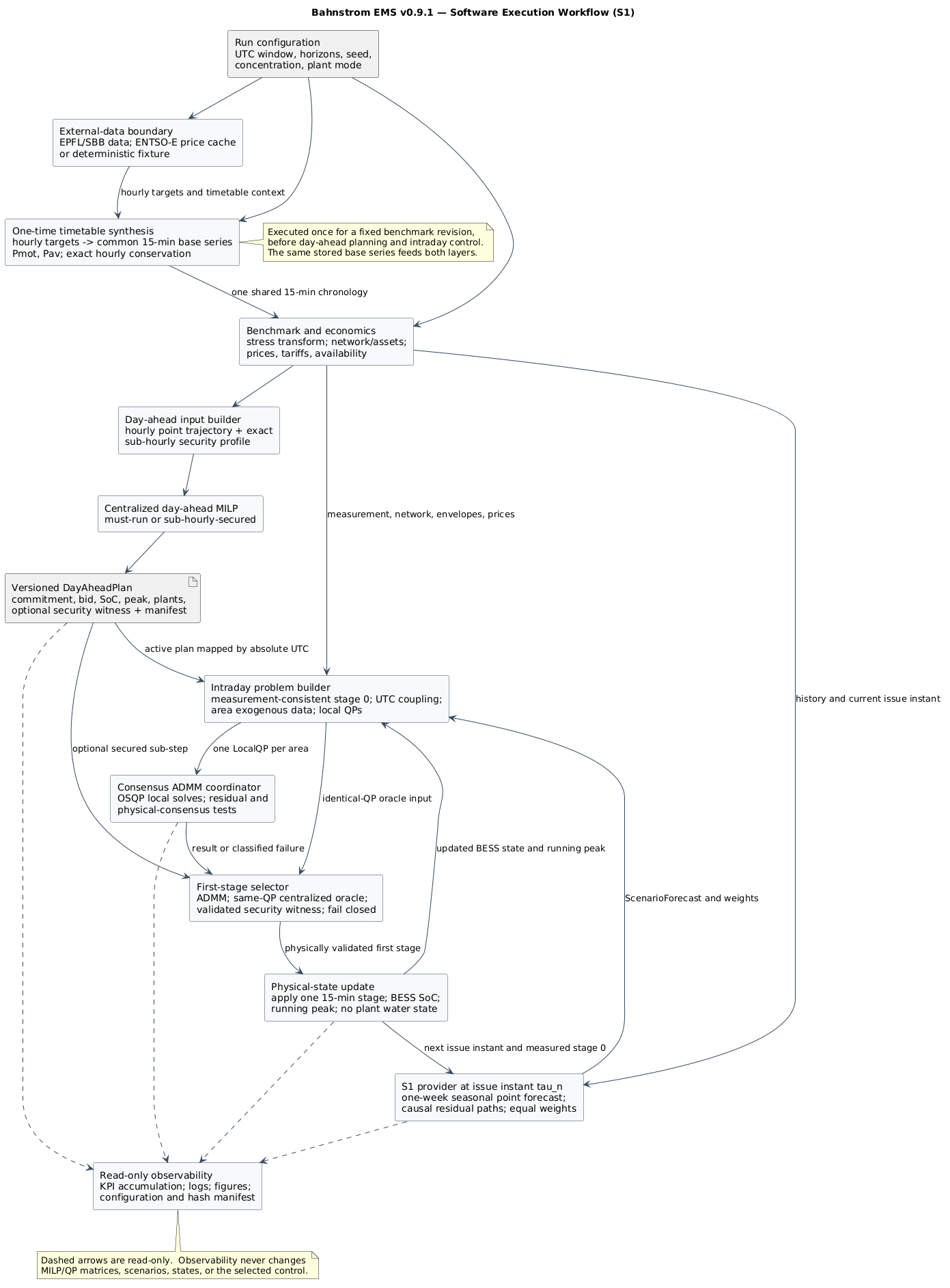}
\caption{Software execution workflow of \software{} v0.9.1 for the S1 configuration.  The
timetable synthesis is executed once before both EMS layers.  At each issue time, S1 is
reissued, stage zero is measured, one action is applied, and the physical state feeds the next
control instant.  Dashed paths are read-only observability.}
\label{fig:software-workflow}
\end{figure}

\Cref{fig:software-workflow} is the software execution workflow for the S1 configuration used
in this paper.

\subsection{Chronological call sequence}

\begin{landscape}
\begin{figure}[H]
\centering
\includegraphics[width=0.98\linewidth,height=0.78\textheight,keepaspectratio]
{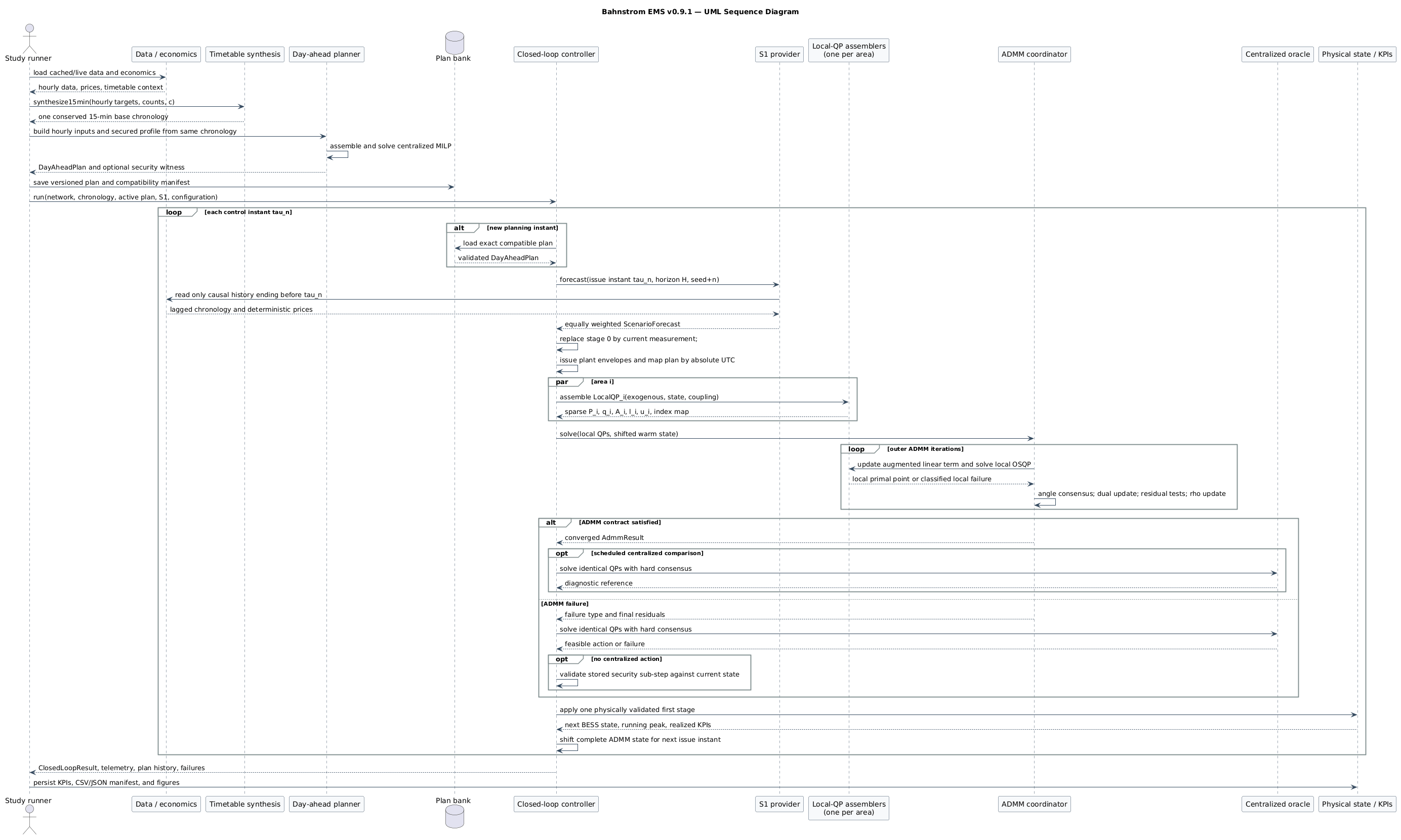}
\caption{UML sequence diagram for one prepared study and its repeated control instants.  The
parallel local-QP call denotes logically independent area assembly and solves.}
\label{fig:sequence-diagram}
\end{figure}
\end{landscape}

The UML sequence diagram in \cref{fig:sequence-diagram} makes the call order and alternative
failure branch explicit.  In particular, the centralized same-QP calculation occurs routinely
only when scheduled for comparison, or exceptionally after an ADMM failure when the relevant
diagnostic/recovery switch is active.  It is not a hidden centralized step inside every
distributed iteration.

\clearpage
\clearpage

\section{Software Architecture and v0.9.1 Module Realization}
\label{sec:architecture}

To get better overview of the software architecture and data modeling, the UML component architecture, a summary of the Python modules and submodules, and the conceptual UML class structure are presented. 

\subsection{UML component architecture}

\begin{landscape}
\begin{figure}[H]
\centering
\includegraphics[width=0.98\linewidth,height=0.78\textheight,keepaspectratio]
{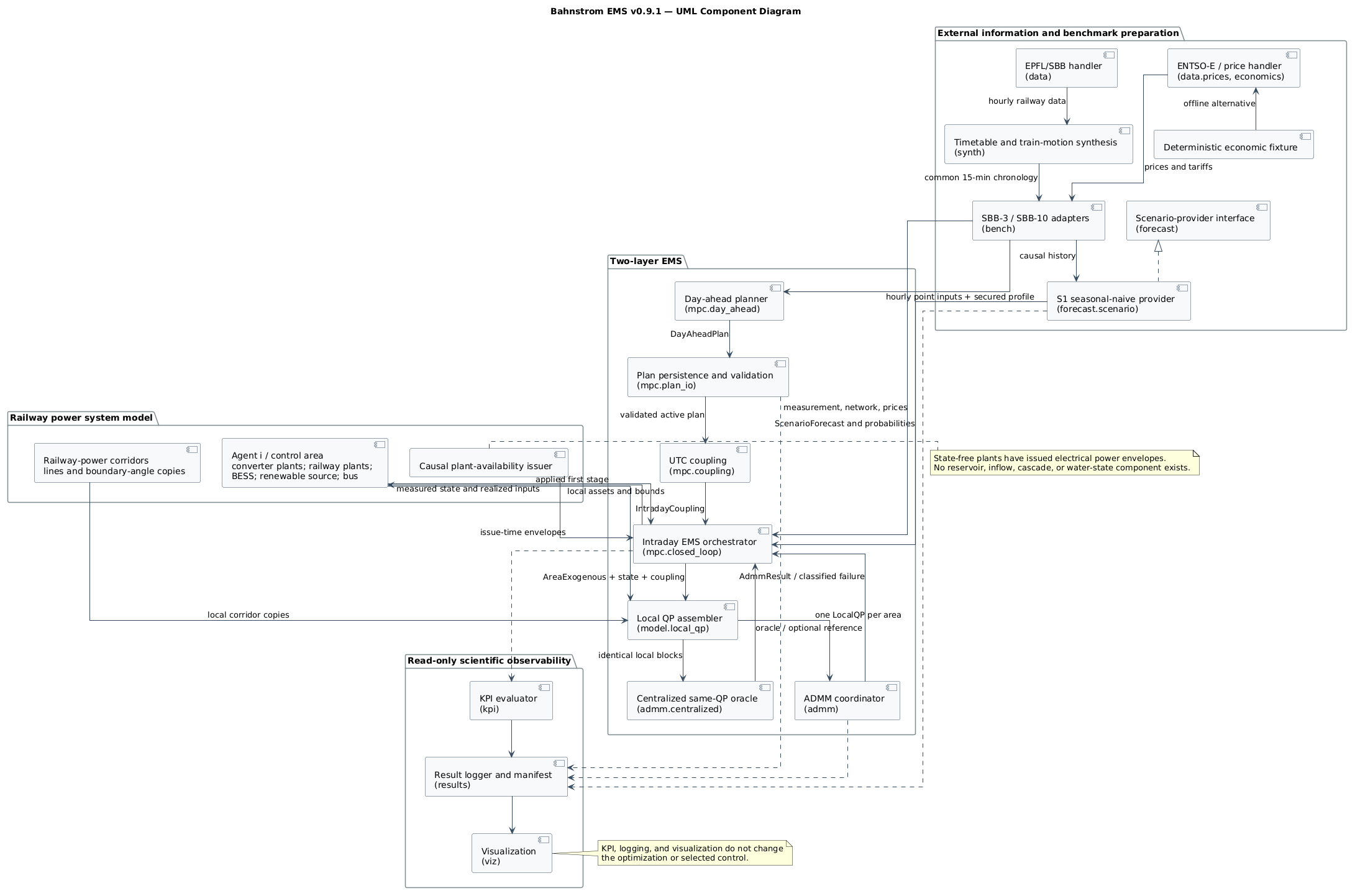}
\caption{UML component architecture of \software{} v0.9.1~\cite{noroozi2026bahnstromems}~\cite{noroozi2026bahnstromems}.  The active forecast component in
Part~II is S1.  An agent is one electrical control area together with its local converter
plants, railway-owned power plants, BESS, renewable source, bus variables, and incident corridor
copies.}
\label{fig:component-diagram}
\end{figure}
\end{landscape}

\Cref{fig:component-diagram} shows the principal executable components and their permitted
dependencies.  The intraday EMS orchestrator coordinates information.
Area assets and corridors enter through the common \texttt{model} objects, while the forecast provider and day-ahead planner cross explicit data
contracts.  The centralized same-QP oracle receives the same local-QP blocks as ADMM.  KPI,
logging, and visualization dependencies point away from the control path.

\clearpage
\clearpage

\subsection{Module and submodule responsibilities}

\Cref{tab:module-architecture} is the source-tree map used for v0.9.1.  It lists the modules
that participate in executable studies, including the power plant increment.  The table
describes ownership rather than call order; call order is already fixed by
\cref{alg:complete-loop,fig:sequence-diagram}.

{\small
\begin{longtable}{@{}p{3.3cm}p{11.6cm}@{}}
\caption{Software organization of \software{} v0.9.1.}
\label{tab:module-architecture}\\
\toprule
Module / submodule & Responsibility \\
\midrule
\endfirsthead
\multicolumn{2}{@{}l}{\small\itshape Table~\thetable\ continued from the previous page}\\
\toprule
Module / submodule & Responsibility \\
\midrule
\endhead
\midrule
\multicolumn{2}{r@{}}{\small\itshape Continued on the next page}\\
\endfoot
\bottomrule
\endlastfoot
\texttt{model}
& Physical and optimization primitives: \texttt{Line}, \texttt{Converter},
\texttt{Bess}, \texttt{RenewableSite}, \texttt{RailwayPlant},
\texttt{PlantAvailability}, \texttt{Area}, \texttt{Network}, and
\texttt{HorizonSpec}.  \texttt{assemble\_local\_qp} maps one area's assets, weighted
scenario data, plant envelopes, loss cuts, and day-ahead coupling into
\texttt{LocalQP}$(P,q,A,l,u,\texttt{idx})$; \texttt{VariableIndex} is the executable
decision-vector map. \\[1mm]

\texttt{admm}
& Builds the owner/copy topology; implements \cref{alg:admm}, warm-state shifting,
classified local-solver errors, residual histories, physical consensus checks, and adaptive
penalty; and supplies the stacked centralized hard-consensus solver used for equivalence and
failure diagnosis. \\[1mm]

\texttt{mpc.day\_ahead}
& Hourly unit-commitment MILP, optional converter must-run map,
optional continuous sub-hourly security block, HiGHS solve, and decoding into
\texttt{DayAheadPlan}.  It owns hourly converter/BESS/plant schedules and the constructive
quarter-hour security witness. \\[1mm]

\texttt{mpc.coupling}
& Absolute-UTC bridge from \texttt{DayAheadPlan} to the present prediction horizon: frozen
commitment, optional gate-closed energy means, BESS references and terminal floor, and peak
target.
Power plant schedules are mapped as audit references and do not become hard
intraday constraints. \\[1mm]

\texttt{mpc.plan\_io}
& Versioned frozen-plan persistence, UTC-index reconstruction, optional security/plant table
persistence, and strict compatibility-manifest validation before replay. \\[1mm]

\texttt{mpc.closed\_loop}
& Two-clock orchestrator of \cref{alg:complete-loop}: plan installation, S1 request and
validation, measured stage-zero overwrite, issue-time plant envelopes, BESS-state injection,
coupling, local-QP assembly, warm-state shift, ADMM, centralized diagnostics, validated
fallback, first-stage application, state/KPI update, and result construction. \\[1mm]

\texttt{mpc.central}
& Separate single-process intraday API retained for development and diagnostics.  It is not
the day-ahead MILP and is distinct from the exact stacked same-QP solver owned by
\texttt{admm}. \\[1mm]

\texttt{synth}
& Davis-based representative train kinematics, category templates, area/bin allocation, and
the timetable-driven, energy-conserving 15-minute construction in \cref{alg:synthesis}. \\[1mm]

\texttt{data}
& External-data boundary for the EPFL/SBB Zenodo records \cite{theiler2026}, cache management,
and ENTSO-E A44 day-ahead-price acquisition and UTC parsing. \\[1mm]

\texttt{economics}
& \texttt{EconomicConfig}, live/fixture selection, signed price transformation, demand-charge
treatment, area-to-zone mapping, and optional loss calibration. \\[1mm]

\texttt{bench}
& SBB-3, SBB-3-15min, and SBB-10 data/network construction; stress transformations; asset and
corridor constants; day-ahead input/security builders; causal plant availability; and
benchmark-specific S1 residual paths.  The plant-free mode is an explicit experimental
configuration, not an implicit missing-data fallback. \\[1mm]

\texttt{forecast}
& Forecast-origin and scenario contracts, validation, S1 provider, and provider interface.
The current repository also contains optional forecast extensions, but they are outside the
active S1 execution workflow and performance claims of this paper. \\[1mm]

\texttt{kpi}
& Per-area and network aggregation of realized market, converter, BESS, renewable,
regenerative, power plant, peak, and solver quantities, including failure-aware ADMM
rates and iteration means. \\[1mm]

\texttt{results}
& Read-only scientific persistence: normalized tick/area, day-ahead, security, scenario, and
plant-action tables; configuration and failure summary; schema information; byte sizes; and
SHA-256 content hashes.  Non-finite JSON values become \texttt{null}; credentials are not
collected. \\[1mm]

\texttt{viz}
& Read-only Matplotlib figures for day-ahead decisions, intraday dispatch, BESS state,
regenerative recovery, plants, scenarios, ADMM/centralized validation, and benchmark topology.
It does not import into the optimization path. \\[1mm]

\texttt{examples}
& Executable study runners and focused demonstrations.  They translate command-line study
parameters into typed configuration, prepare or load plan banks, choose the S1 provider, invoke
the common closed loop, and request result/figure formats. \\
\end{longtable}
}

\subsection{Conceptual UML class structure}

\begin{landscape}
\begin{figure}[H]
\centering
\includegraphics[width=0.98\linewidth,height=0.78\textheight,keepaspectratio]
{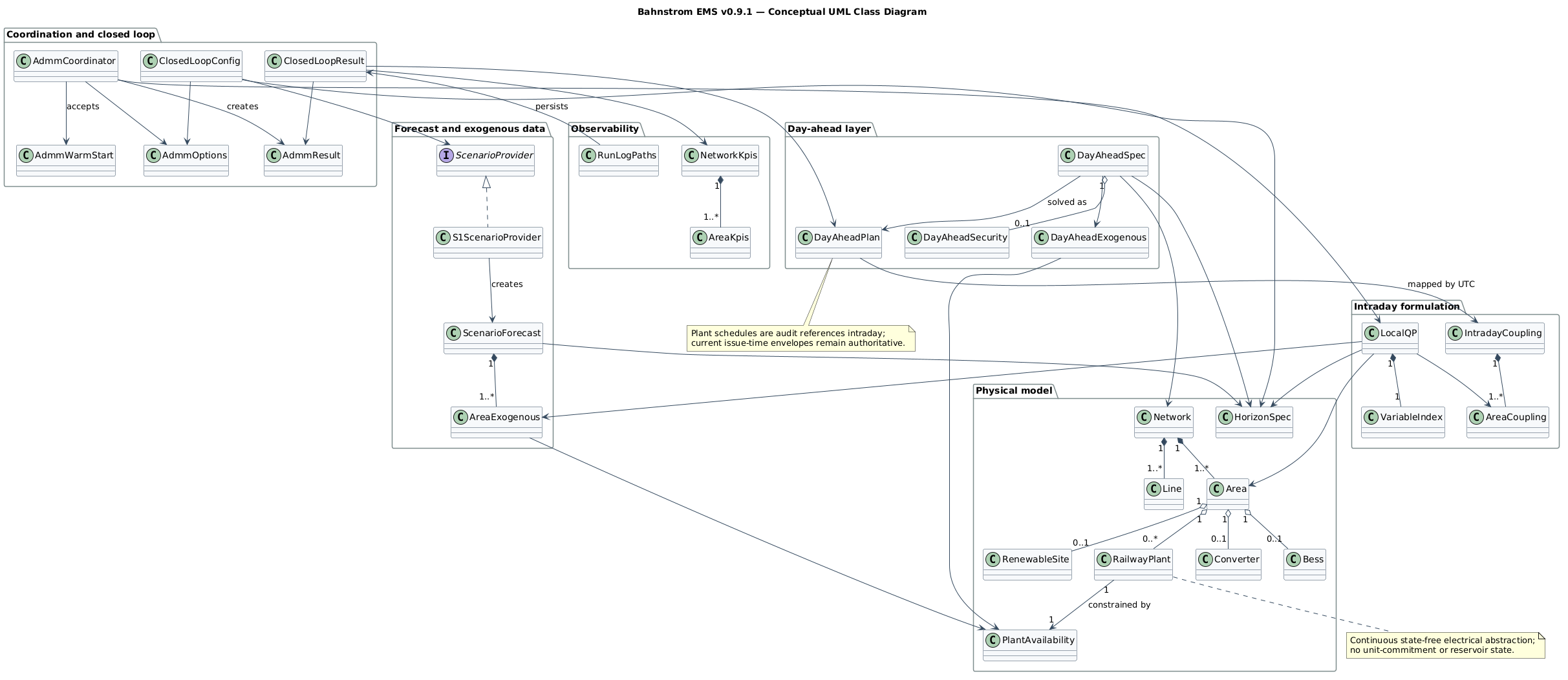}
\caption{Conceptual UML class diagram for \software{} v0.9.1.  Only the key physical,
forecast, planning, QP, coordination, closed-loop, and observability classes are shown.}
\label{fig:class-diagram}
\end{figure}
\end{landscape}

The class diagram in \cref{fig:class-diagram} illustrates the conceptual design.
Composition diamonds identify owned physical or result structure;
ordinary dependencies identify data flow.  The \texttt{ScenarioProvider} interface decouples
forecast formation from the controller, while \texttt{LocalQP} and \texttt{VariableIndex}
decouple mathematical assembly from ADMM numerical state.  

\clearpage
\clearpage

The central chain of ownership is as follows.  A \texttt{Network} composes areas and lines; an
\texttt{Area} aggregates only the assets it physically contains.  A forecast provider creates a
\texttt{ScenarioForecast}, which composes area exogenous blocks.  A
\texttt{DayAheadPlan} is mapped to an \texttt{IntradayCoupling}, which composes area couplings.
Each \texttt{LocalQP} refers to one area, one area exogenous block and coupling, and owns its
variable-index map.  The coordinator accepts local QPs and an optional warm state and creates
an \texttt{AdmmResult}.  The closed-loop result retains the applied outcome, plan and numerical
records from which network and area KPIs are derived and persisted.

\begin{singlespace}
\bibliographystyle{IEEEtran}
\bibliography{Bahnstrom_EMS_P1_Common}
\end{singlespace}

\end{document}